\documentclass{article}
\usepackage[utf8]{inputenc}
\usepackage[T1]{fontenc}
\usepackage{amsmath}
\usepackage{amsthm}
\usepackage{amssymb}
\usepackage{tikz-cd}
\usepackage{authblk}

\usepackage{tikz}
\usepackage{subcaption}
\usepackage{graphicx}
\usetikzlibrary{decorations.pathmorphing}
\usepackage{float}
\usepackage{cancel}
\usepackage{listings}
\usepackage{cite}
\usepackage[colorlinks=true, linkcolor=blue, urlcolor=blue]{hyperref}
\usepackage{footnote}
\usepackage{nicefrac}

\theoremstyle{definition}

\theoremstyle{remark}

\numberwithin{equation}{section}

\usepackage{tablefootnote}

\begin{document}
\title{Why Do Students (Not) Choose Second-Cycle Mathematics Studies? Questionnaire and Graduate-Tracking Evidence from Poland}
\author[1,4]{Filip Turoboś\thanks{filip.turobos@p.lodz.pl}}
\author[2]{Żywilla Fechner\thanks{zywilla.fechner@p.lodz.pl}}
\author[2]{Jacek Stańdo\thanks{jacek.stando@p.lodz.pl}}
\author[3]{Nicole Meisner\thanks{meisner.nicole@gmail.com}}
\affil[1]{Institute of Mathematics, Lodz University of Technology}
\affil[2]{Center of Mathematics and Physics, Lodz University of Technology}
\affil[3]{Faculty of Mathematics, Informatics and Mechanics, University of Warsaw}
\affil[4]{Institute of Theoretical and Applied Informatics, Polish Academy of Sciences}
\date{May $4^{\text{th}}$, $2026$}

\maketitle


\begin{abstract}
The Bologna Process has substantially reshaped higher education systems across Europe, including the structure and organization of mathematical studies in Poland. One of the increasingly visible consequences of these transformations is the relatively low retention rate between first- and second-cycle studies. The aim of this paper is to investigate selected factors associated with students' willingness to continue mathematical education after obtaining a Bachelor's degree.

The study combines questionnaire-based research conducted among mathematics students from 13 Polish higher education institutions with an auxiliary analysis of nationwide graduate-tracking data obtained from the Polish Graduate Tracking System (ELA). The survey investigated students' opinions on general and specialized courses, perceived labour-market usefulness of studies, organizational aspects of university education, specialization satisfaction and future educational intentions.

The results indicate that students willing to continue second-cycle studies generally evaluate both the substantive quality and practical usefulness of their studies more positively than students intending to leave mathematics or change institutions. Satisfaction with the chosen specialization emerged as one of the strongest differentiating factors between the analysed groups. At the same time, a substantial proportion of respondents expressed doubts regarding the professional utility of continuing mathematical education, despite administrative labour-market data suggesting several advantages associated with obtaining a Master's degree. The findings suggest that retention in mathematics is shaped not only by academic difficulty, but also by the perceived relationship between university curricula, specialization structures and labour-market expectations. We conclude with several recommendations regarding curriculum design, cooperation with external stakeholders and other aspects of second-cycle mathematical education.
\end{abstract}

\section{Introduction}

On 2024, there was 31 universities in Poland which offered the access to the first-cycle mathematical studies, out of which 27 provided students with second-cycle mathematical studies as well. Despite this rather broad selection, the second-cycle studies are often overlooked by the students of mathematics, as well as the students of other subjects, as pointed out in data collected by OPI-PIB ELA.\footnote{National Information Processing Institute --
National Research Institute, offering access to aggregated data obtained by Polish Graduate Tracking System, which monitors the professional careers of polish students.
\href{https://www.ela.nauka.gov.pl/en}{https://www.ela.nauka.gov.pl/en}} This is not a new problem -- in general, there is a lot of papers on drop-out ratio in STEM-related fields \cite{burke2019student,tight2020student,valencia2024dropout}, but unclear definition of a drop-out \cite{xavier2020literature} significantly hampers most literature review attempts in this direction. 

To ramp up the difficulty of retention problems in mathematics in particular, a large part of the community of academic teachers have noticed that ``\textit{Mathematics students are no longer what they used to be.}'' -- they tend to avoid more abstract courses in favour of more data-analysis- or actuarial-oriented subjects.

There were several reports on the topic of retention between the first and the second cycle of studies (e.g. \cite{rotem2021dropping}). In terms of retention, there are papers which are dedicated specifically to retention in mathematics (cf. \cite{fenwick2009recruitment,Geisler2018}), either as a part of the curriculum or as a study field. Individual aspects are investigated as well, usually in the context of STEM-related fields (see, e.g. \cite{Lytle2023STEMEngagement}). There also exists a fair number of mathematical models which allow to predict the student retention (see \cite{fenwick2009recruitment,kilian2020predicting,burke2019student}), 
none of these took into account all of the following factors at once:
\begin{itemize}
\item pinpoints the problem of transition between the Bachelors' and Masters' degree studies;
\item focusing solely or mainly on mathematics;
\item providing reasons for that situation which would be backed up by data gathered from students.
\end{itemize}

The main aim of this paper is to fill that gap in the following ways: we aim to identify the main reasons for the students to finish only the first-cycle studies and simultaneously explore the potential methods of improving the percentage of Bachelors who stay at the university for Masters degree. We try to explore the nature of expectations that upcoming first-cycle studies graduates have towards the postgraduate courses offered by the universities by conducting a survey amongst them (with the help of HEI authorities).

The paper is organized as follows -- in the following subsection we discuss the main assumptions and implementation of Bologna process with an emphasis on the Polish higher education experience. The sections 2. and 3. discuss the survey methodology and preliminary survey design, respectively. Then, we follow with the description of a study group and analysis of the responses obtained in the questionnaire. This part is divided into several subsections, dedicated to each part of the survey separately, followed by some joint insights. The penultimate main part of this article describes the obtained results and conclusions drawn from the questionnaire outcomes. We also give some recommendations on steps worth undertaking to improve the current situation with second-cycle studies. Lastly, we have funding informations, acknowledgements and references. 

\subsection{On Bologna process in Poland}

While the beginnings of the Bologna process can be dated back to 1998 \cite{SorbonneDeclaration1998}, the reformative process of European higher education has truly begun a year later. This event was marked by signing of the Bologna Declaration \cite{BolognaDeclaration1999} by 29 representatives of ministries of higher education from multiple countries across the Europe. The progress on the implementation of those recommendations along with the new statements/guidelines are discussed on biennial conferences of higher education ministers of participating countries.

The main assumptions of the process were as follows:

\begin{itemize}
\item introduction of a system of transparent and easily comparable degrees by introducing the Diploma Supplement;
\item agreement on degree structure based on three cycles -- bachelor studies, which last between 3 and 4 years, master degree studies lasting a year and a half or two years, followed by the doctorate;
\item promulgation of the European Credit Transfer System;
\item promotion of mobility amongst people related to academia (students, scientists, academic teachers as well as administrative workers);
\item cooperation on the improvement of higher education quality across the Europe;
\item promotion of European Higher Education Area and the interdisciplinary research.
\end{itemize}

The most important milestone for Bologna process in Poland was the introduction of \textbf{The Act of 27 July 2005 on Higher Education} \cite{PrawoOSzkolnictwieWyzszym2005},
 which started the changes towards unification of polish higher education with European standards. It is worth pointing out that the start of the Bologna process in Poland occurred prior to the law change -- for example, the diploma supplements were handed out by universities since the 1st January 2005 \cite{OlenczukPaszel2009}.

It is hard to evaluate the impact of the Bologna process solely on the didactics of mathematics, due to the wide scope of its effect, which incorporates changes in the quality, universality of the education as well as the social dimension of studying \cite{Crosier2013,Usaci2012,Huisman2012,Terry2008}. This provides us (as well as other researchers) with a neat unexplored area of how did the mathematical education has evolved over the course of this process.

\section{Methodology}

A cross-sectional research design including a mix of quantitative and qualitative descriptive approaches were used. The initial questionnaire was devised based on the questions most commonly asked in other similar studies, while simultaneously utilizing some suggestions from the discussions with the students of Lodz University of Technology. The preliminary survey with open-ended questions was conducted amongst people who have finished mathematics at Lodz University of Technology with either Masters or just Bachelor degrees, to check for the most recurring patterns in the answers. This is one of the well-known techniques for on-line survey pre-testing \cite{kaczmirek2017higher,geisen2020compendium}, which is often underestimated technique with immense impact on the quality of the final survey -- see \cite{hilton2017importance,kaplan2020multi,neuert2021effects,rothgeb2007questionnaire}, as well as references therein. This pilot version of survey was conducted over May of 2024.

After a careful analysis of the pilot survey results, the close-ended questions were prepared and the final version of  questionnaire was then sent to university authorities of every HEI in Poland which offered courses on mathematics, with a request of distributing them amongst the students of second and third year of mathematical studies in Poland. In around $40\%$ of cases we have received a positive response from the universities. The responses were collected via online survey with the access limited to the owners of e-mails within the \textsl{edu} domain. This (along with the distribution of the questionnaire solely via universities inner communication channels) enabled us to partially protect ourselves from the situations, where the respondent did not belong to the target group. To increase the response rate, the survey was designed to take about 10 minutes to complete. The questions of more personal nature (such as household income, GPA etc.) were left non-obligatory. This resulted in us being able to collect the data from 170 students, scattered across 13 universities. The collection of results took place between June and October of 2024.

\begin{figure}
\centering
\includegraphics[width = 0.55\textwidth]{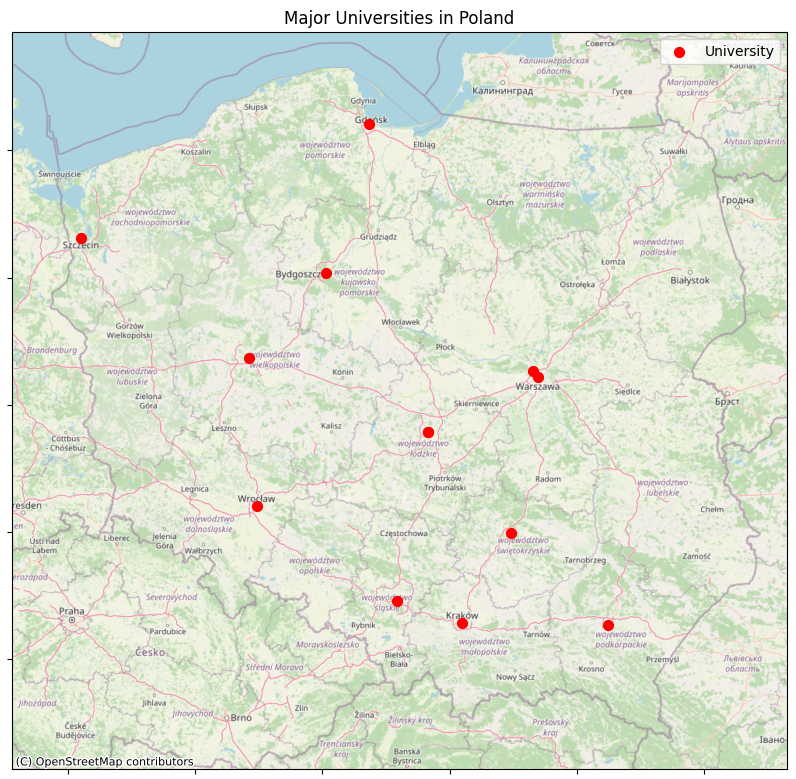}
\caption{Approximate positions of the universities whose students participated in the survey.}
\end{figure}

Lastly, the data analysis process was conducted on the answers dataset. We utilized Python 3.12.5 scientific environment with Numpy \cite{harris2020array}, Pandas \cite{pandas}, Scipy \cite{scipy}, Matplotlib \cite{matplotlib}, Jupyter \cite{jupyter}, Seaborn \cite{Waskom2021} and Plotly \cite{plotly} libraries.

Descriptive statistics were used to summarize the contents of the responses and obtain insights on the predominant reasons behind the decision whether to embrace Masters degree studies. We have also applied correlation analysis and statistical testing for exploring the relationships between the opinions on particular subjects. 
Additional data-mining algorithms (such as Apriori algorithm -- see \cite{du2016research} as well as the original paper introducing this method \cite{agrawal1993mining}) were utilized to discover additional associations and identify frequently coupled choices both in singular, multiple-choice questions as well as the answers across the several different questions.

\section{Preliminary survey design}

Below we describe the questions included in the preliminary survey. Open-ended questions used to gather exemplary answers were marked with boldface font. In the following we omit the questions related to formal agreements for personal information processing etc.

\begin{table}[H]
\begin{tabular}{ |p{0.5\textwidth}|p{0.5\textwidth}| } 
 \hline
 \large{Native formulation} & \large{English formulation} \\\hline\hline
 Jak oceniasz swoje dzisiejsze samopoczucie? & How are you feeling today?\\\hline
Uczelnia na której studiuję to: & The university I study/studied on is: \\\hline
\textbf{Głównymi czynnikami, którymi kierowałem się przy wyborze uczelni były:} & \textbf{The main factors you considered when choosing your university were:}\\\hline
\textbf{Co miało wpływa na Twoją decyzję o podjęciu studiów matematycznych?} & \textbf{What influenced your decision to pursue the field of mathematics?}\\\hline
\textbf{Jaka jest Twoja opinia na temat przedmiotów ogólnych, których uczysz się na studiach?} & \textbf{What is your opinion on general, obligatory courses during your studies?}\\\hline
\textbf{Jaką specjalność wybrałeś/-aś w ramach toku studiów?} & \textbf{What specialization did you choose during the first cycle of your studies?}\\\hline
\textbf{Jakie powody stały za tym wyborem specjalności/przedmiotów obieralnych?} & \textbf{What influenced your decision to choose this specialization/elective courses?}\\\hline
\textbf{Co uważasz na temat wybranych przez siebie przedmiotów specjalnościowych i obieralnych?} & \textbf{What is your opinion on specialized/elective courses you chose during your studies?}\\\hline
\textbf{Jaka jest Twoja ocena organizacyjnej strony Twoich studiów?} & \textbf{What is your opinion on the organizational side of the university (schedules, compatibility with part-time work, flow of information etc.)?} \\\hline 
Jeżeli mógłbyś/mogłabyś zmienić coś od strony organizacyjnej ze studiami na swojej uczelni, to co by to było? & If you were able to change one thing regarding the organization of your university -- what would that be? \\\hline
\textbf{Czy masz zamiar kontynuacji studiów na drugim stopniu na Twojej obecnej uczelni? Jeśli nie -- jaka jest planowana przez Ciebie dalsza ścieżka kariery?} & \textbf{Do you plan to stay at the same university for the second-cycle of studies? If not, what is your alternative choice?} \\\hline
\textbf{Jakie czynniki mogły mieć wpływ na Twoją decyzję w tej sprawie?} & \textbf{What are your motivations behind this decision?} \\\hline
Czy którekolwiek z pytań było według Ciebie nieprecyzyjnie sformułowane lub niejasne? & Was any of the questions formulated in imprecise or confusing way?\footnotemark\\\hline
\end{tabular}
\caption{The questions included in the pretesting questionnaire.}
\end{table}
\footnotetext{Due to the probing nature of this question, it was not included in the subsequent, proper survey.}

The open-ended nature of each question allowed us to formulate most often selected choices for each of the questions. Due to sheer amount of creativity expressed in the pilot version of this survey, we have decided that the question regarding the \textit{one improvement to the organizational side of the university} should remain open.

\section{Questionnaire results and their analysis}

We begin with a brief information about the demographic structure of the analysed sample. Since these questions were not obligatory, only $90$ participants of the survey provided answers to these. Nevertheless, below we provide a brief, visual summary of collected data. In order to anonymize the data, the individual information on participants will not be included in publicly available dataset.

\begin{figure}[H]
    \centering
    \begin{subfigure}{0.45\textwidth}
        \includegraphics[width=\linewidth]{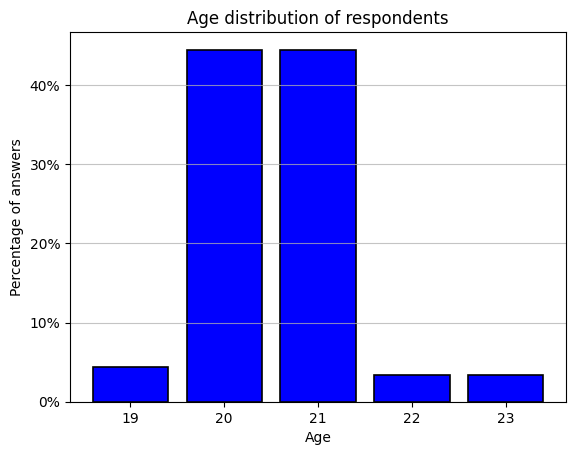}
        \caption{Distribution of age amongst the respondents (mean: 20.56, standard deviation: 0.78).}
        \label{fig:image1}
    \end{subfigure}
    \hfill
    \begin{subfigure}{0.45\textwidth}\vspace{-0.5cm}
        \includegraphics[width=\linewidth]{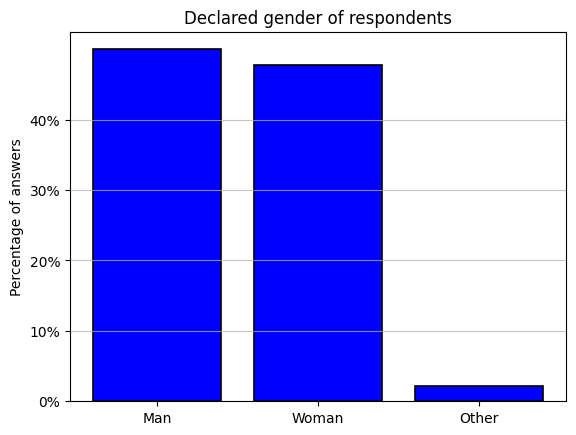}
        \caption{Declared gender of the respondents}
        \label{fig:image2}
    \end{subfigure}
    
    \vspace{0.5cm} 
    
    \begin{subfigure}{0.45\textwidth}
        \includegraphics[width=\linewidth]{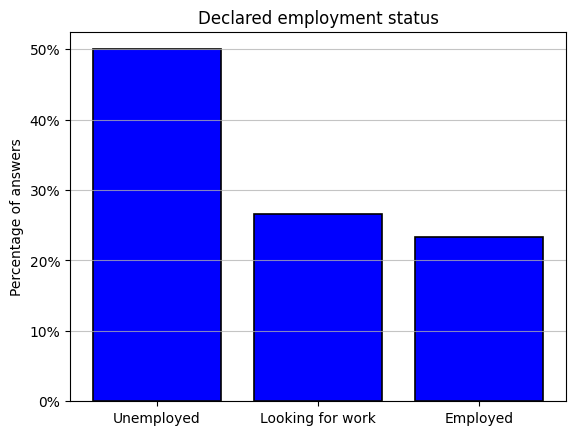}
        \caption{Employment status (self-reported).}
        \label{fig:image3}
    \end{subfigure}
    \hfill
    \begin{subfigure}{0.45\textwidth}
        \includegraphics[width=\linewidth]{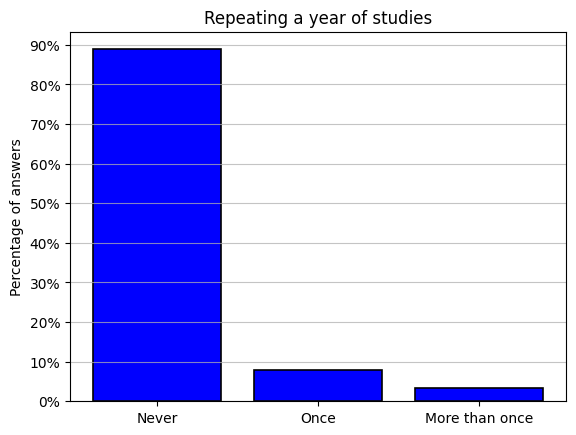}
        \caption{Declared number of retaken study years.}
        \label{fig:image4}
    \end{subfigure}
    
    \caption{Basic parameters of the sample.}
    \label{fig:grid}
\end{figure}

A few key things can be said about the presented data. Almost $50\%$ of respondents are not employed, nor they are interested in undertaking a permanent job in the nearest foreseeable future. This is mostly the case of students, who plan on staying on the university for the second-cycle studies, as we will see in the latter part of the paper. The other half of respondents either had their first-hand experience with the job market, or they will be having it very soon. 

As far as the working students are concerned, on average they work $20.05$ hours per week, which corresponds to $\nicefrac{1}{2}$ of standard weekly working time in Poland. Predictably, this quantity bears a rather large standard deviation, which equals $12.84$. For more details, please refer to the following histogram:

    \begin{figure}[H]
    \centering
        \includegraphics[width=0.5\linewidth]{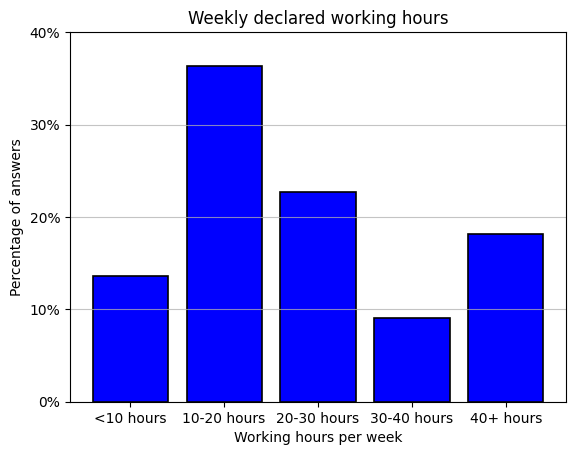}
        \caption{Declared number of weekly working hours.}
        \label{fig:image123}
    \end{figure}
    
As far as the financial situation of the students is concerned, most of them do not have to worry about finances. We asked the students to rate their need to control their own spendings in relationship to the money they have using 10 level scale. The lowest value 1 corresponded to the situation where almost every expense has to be planned and 10 was used to describe the case, where no saving behaviour was needed to continue the current lifestyle.

\begin{figure}[H]
\centering
\includegraphics[width=0.5\linewidth]{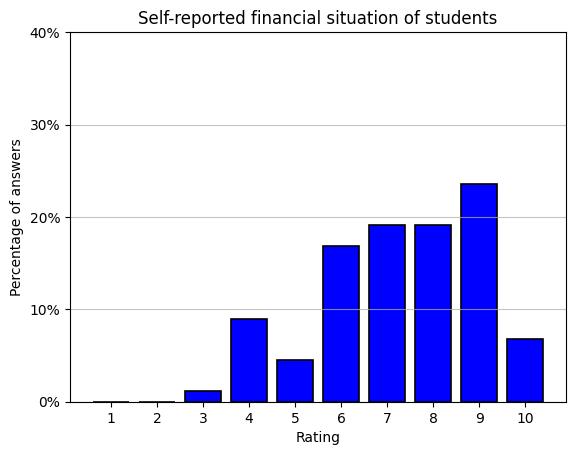}
\caption{Self-reported financial status -- mean rating: 6.30, standard deviation: 1.75.}
\end{figure}

As far as the living conditions of students are concerned, most of them study near their parental home or flat. Only scant few own a living place on their own or with their spouse. Those who are not lucky enough to have their living place near the university tend to rent apartments on their own. Access to dormitories is often limited, but some students avoid these due to the cultural and societal aspects related to staying at these (necessity of having a room-mate, parties taking place therein etc.).\footnote{Those observations were not taken directly for the collected data, they are based on a few conversations with individual students, who were curious about this research.}

\begin{figure}[H]
\centering
\includegraphics[width=0.7\linewidth]{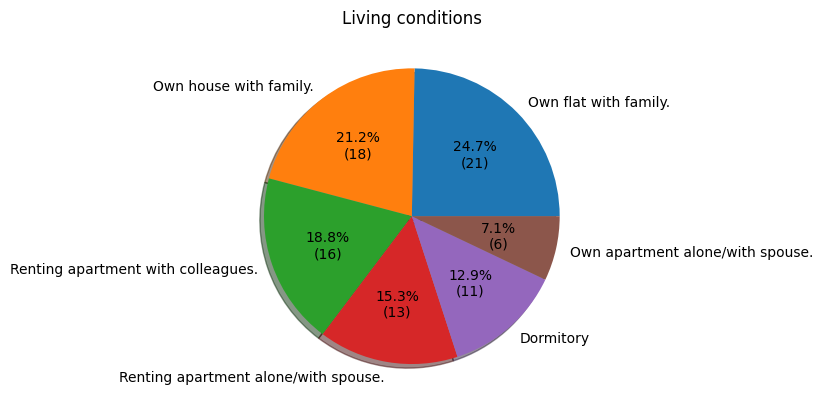}
\caption{Self-reported living conditions.}
\end{figure}

The respondents GPAs were on average very close to 4.0, with standard deviation of 0.51.\footnote{In Poland, higher education operates on grades from 2 to 5, with the inclusion of 3.5 and 4.5. Grades below 3 are not considered being a passing grades, 5 is the best attainable grade.}

\begin{figure}[H]
\centering
\includegraphics[width=0.5\linewidth]{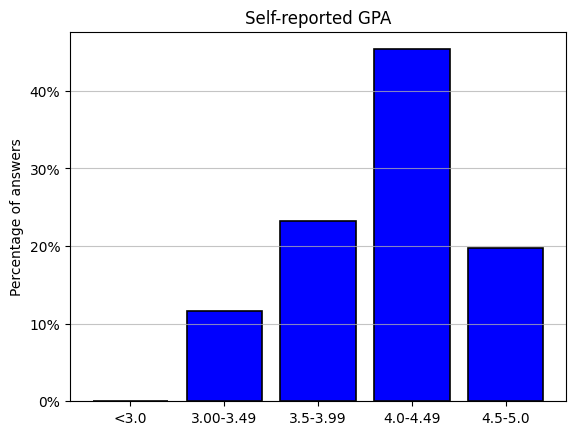}
\caption{Self-reported average from grades obtained in the last semester.}
\end{figure}

On average, the self-reported financial standings of students who live with their families in houses tended to be higher than the rest of the group, though the difference was not significant in most cases.
    
Lastly, unlike most of the strictly engineering and science-oriented fields (see \cite{hill2010so,silbey2016so}) there seems to be similar number of women and men who are attracted to studying mathematics. In the Subsection \ref{subsec:genderr} we will explore to what extent does the gender correlate with the opinion on particular statements regarding studying.

\subsection{Opinion on general courses}

Before discussing any of the findings for this section, let us present the barplot depicting the distribution of students' answers to each of the questions in this part of the questionnaire.

\begin{figure}[H]
\includegraphics[width = \textwidth]{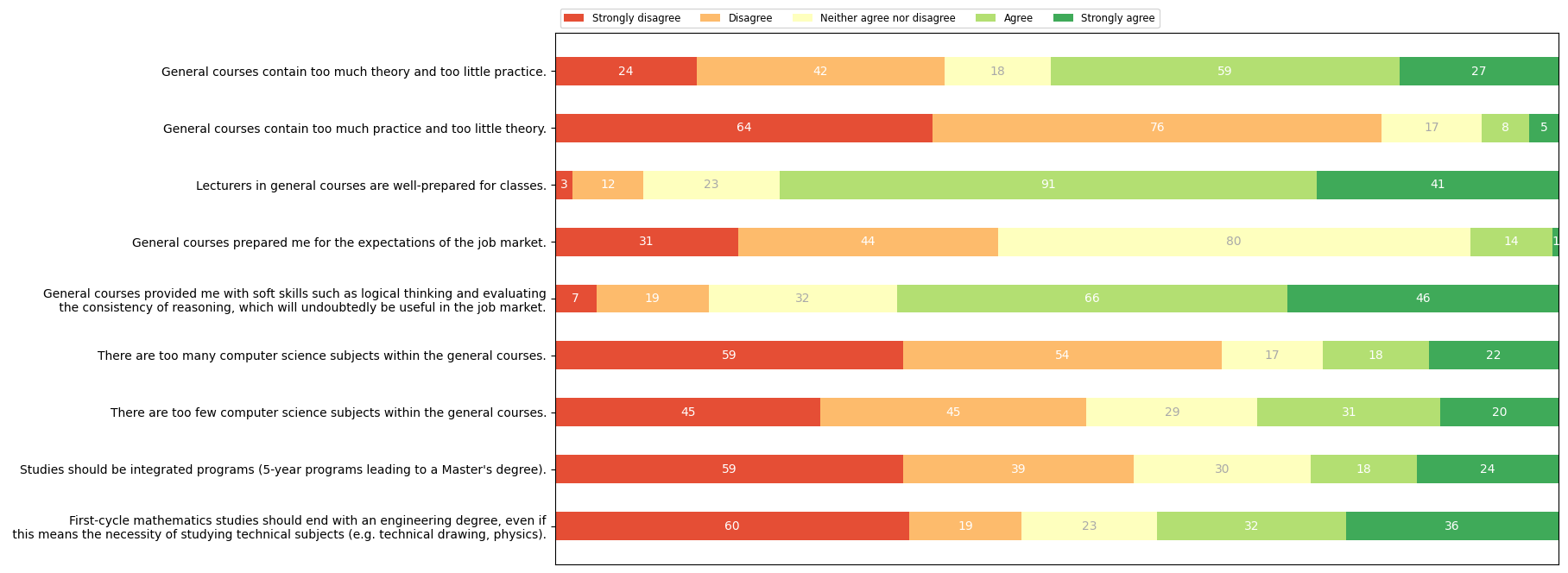}
\caption{The distribution of students' opinions on general courses.}
\end{figure}

For the purpose of statistical analysis, the answers were translated from standard 5-degree Likert scale\footnote{Consisting of answers "Strongly disagree", "Disagree", "Neither agree nor disagree", "Agree", "Strongly agree".} to the values between 1 and 5, where 1 corresponded to strong disagreement and 5 was mapped to the opposite end of the scale. The following table contains a summary on average opinion on each subject and the standard deviation (as a measure of students' divergences on views on a particular matter).

\begin{table}[H]
\begin{tabular}{ |p{0.8\textwidth}|p{0.1\textwidth}|p{0.1\textwidth}| }
\hline
\textbf{Question} & Mean & SD \\\hline
General courses contain too much theory and too little practice. & 3.135 & 1.337 \\\hline
General courses contain too much practice and too little theory. & 1.906 & 0.962 \\\hline
Lecturers in general courses are well-prepared for classes. & 3.912 & 0.903 \\\hline
General courses prepared me for the expectations of the job market. & 2.471 & 0.905 \\\hline
General courses provided me with soft skills such as logical thinking and evaluating
 the consistency of reasoning, which will undoubtedly be useful in the job market. & 3.735 & 1.102 \\\hline
There are too many computer science subjects within the general courses. & 2.353 & 1.386 \\\hline
There are too few computer science subjects within the general courses. & 2.624 & 1.359 \\\hline
Studies should be integrated programs (5-year programs leading to a Master's degree). & 2.465 & 1.419 \\\hline
First-cycle mathematics studies should end with an engineering degree, even if
 this means the necessity of studying technical subjects (e.g. technical drawing, physics). & 2.794 & 1.591 \\\hline
\end{tabular}
\caption{Means and standard deviations of responses to questions regarding general courses.}
\end{table}

The findings in this section can be summarized to the following highlights:
\begin{itemize}
\item Students definitely agree on the fact, that their lecturers are well prepared for the classes (over $77\%$ of students agree on that matter);
\item They are also highly convinced that the soft skills obtained via the general courses will be of great use in the job market (over $65\%$ of students either agree or strongly agree with this statement);
\item At the same time they doubt whether general courses provide them with skills which are expected of them at the job market. Over $90\%$ have selected non-positive responses for this question, with $47\%$ of undecided respondents;
\item There is a noticeable fraction of students who believe that general courses should contain a bit more practical approaches.
\item The majority of students seems to be satisfied with the number of computer-oriented subjects in their curricula, with a small percentage of students who believe the number of such subjects should be increased (exactly $30\%$).
\item Only one in four students believe there might be some merit in returning to the 5-year study program leading directly to the Master's degree.
\item The question which sparks the most divided viewpoints is about the engineering nature of mathematics courses. Almost half of the students (slightly below $47\%$) do not believe that mathematics higher education should award the engineering degree if it obliges them to learn physics and other technical subjects. At the same time, $40\%$ students have positive view on that matter.	
\end{itemize}

The relationship between these opinions and other answers of the respondents will be discussed in the latter parts of the paper.

\subsection{Opinion on specialized courses}

On most universities, the students' are asked to choose a specialization near the start of their fourth semester (although in some cases this takes place earlier). Those HEI which do not offer a selection of specializations in the form of a \textit{labeled package} usually tend to have some form of elective courses nevertheless -- so we have opted for asking them several questions about those \textit{semi-specialized} courses of their choice. In the questionnaire, we have allowed the students to either select one of the possible descriptions of their specialization of choice, but we have also integrated a possibility of providing own description, should such a need arise. The latter option was sparsely used (12 cases), and all of these answers were assigned to the closest descriptions from the set of predefined responses. The specialities of choice for students can be seen in the following chart:

\begin{figure}[H]
\includegraphics[width = \textwidth]{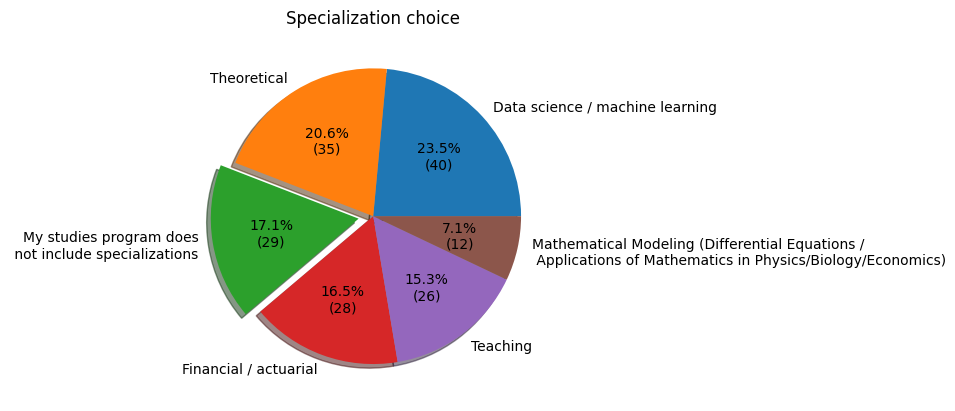}
\caption{The selection of specializations amongst respondents.}
\end{figure}

We can see that most preferred choices (and perhaps the most commonly offered specializations) include theoretical one and specialization in data science. Teaching and actuarial-oriented specializations contribute almost in equal proportion to nearly $40\%$ of choices. The remaining $8.5\%$ of students choose specialities related to mathematical modeling and differential equations. This can be explained to some extent with the fact, that the remaining specialities have rather clearly defined job counterparts and the mathematical modeling might seem either mysterious or simply less appealing than the remaining choices.

\begin{figure}[H]
\centering
\includegraphics[width = 0.85\textwidth]{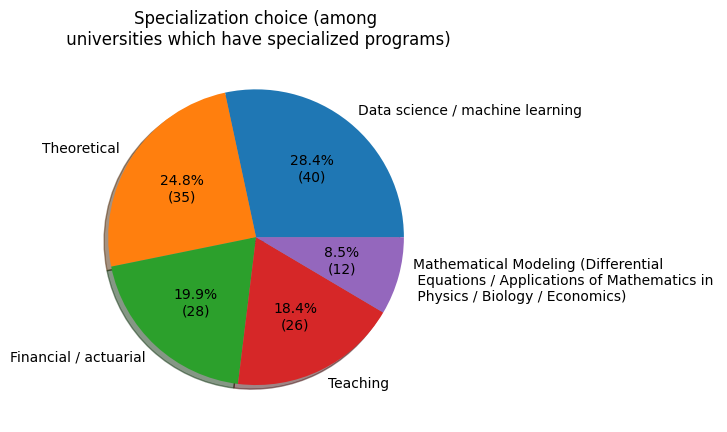}
\caption{The selection of specializations amongst respondents after excluding universities without specialized programs.}
\end{figure}

Proceeding with the brief analysis of the responses to the questions concerning the specialized courses, let us begin with the following visualization:

\begin{figure}[H]
\centering
\hspace{-0.1\textwidth}\includegraphics[width = 1.1\textwidth]{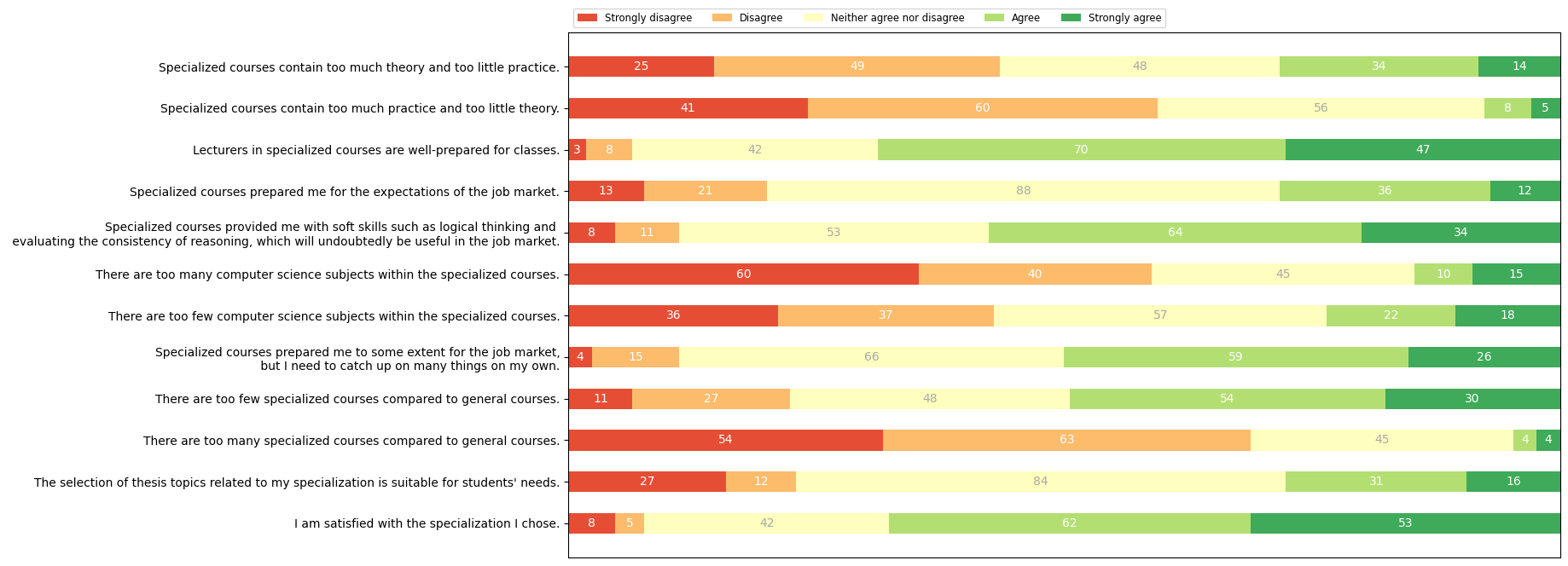}
\caption{The distribution of students' opinions on specialized courses.}
\end{figure}

The good news is that the students tend to respond positively towards the question regarding their satisfaction from the chosen speciality, which means that the fields of study they choose to focus on in most cases meets their expectations. They also tend to appreciate the preparation level of the lecturers for the specialized courses. 

Regarding the balance between the theory and practice, there is a small but noticeable shift towards the opinion that these subjects have sufficient amount of practical component. This is reflected by the increasing number of respondents (as when comparing with the opinions on general subjects) who neither agree nor disagree on the fact, that the evenness between these two constituents is shifted to either theory or practice.

Perhaps the most noticeable change in opinions when comparing the views of students on general and specialized subjects is their usefulness on job  market. The large majority of respondents is convinced that the knowledge and skills they obtain during the speciality courses are either sufficient on their own to meet the expectations of their future employer (almost $30\%$ of responses) or they believe that these constitute to a decent foundation, upon which they will be able to expand the repertoire of their own abilities (half of the questioned students). It is important to point out that over $50\%$ of respondents are unsure whether they will be able to match expectations of the job market with the skills they gain via elective/specialized courses -- this can indicate that they do not know, what the tentative employer might expect from them. 

The results of the questionnaire suggest that the equilibrium between the number of general and elective choices in the case of most universities could be shifted a bit towards increased number of specialized courses. This opinion is shared -- to some extent -- by nearly $50\%$ of the respondents. This is further backed up by the fact that only $22\%$ of students believing there is sufficient or excessive number of specialized courses in their curriculum. The discourse about the number of computer-related subjects amongst the specialized courses tends in similar direction, however there is a stronger opposing fraction towards the increasing the number of such subjects in the elective part of their curriculum. 

\subsection{Retention between first and second-cycle studies}\label{subsec:retentionbetween}

We now turn our attention to the main problem discussed within the paper and that is the retention between the first-cycle studies and Masters studies. 

\begin{figure}[H]
\centering
\includegraphics[width = 0.85\textwidth]{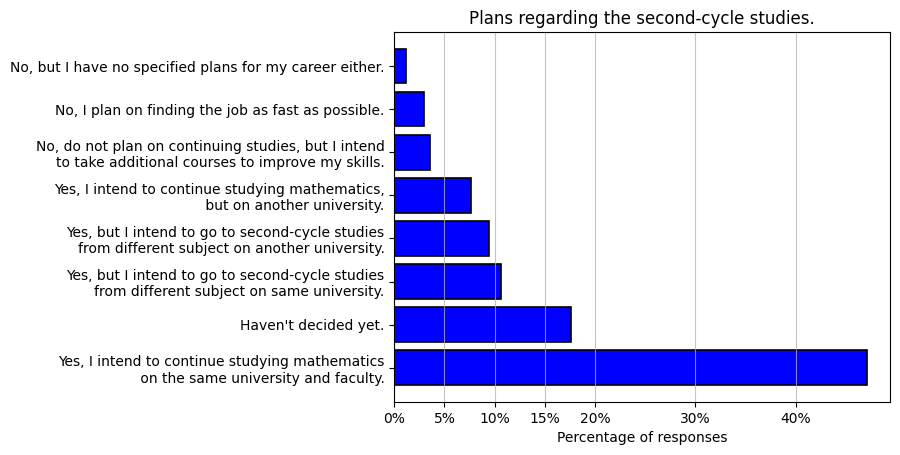}
\caption{Self-reported plans for subsequent career after first-cycle studies.}
\end{figure}

Less than $55\%$ of students intend to continue their mathematical education (with $7.6\%$ indicating, that they are willing to move to another HEI). What is interesting is the fact that a noticeable fraction of students ($20\%$) intend to continue studies, but they want to focus on other topics than mathematics alone. Around $17.5\%$ of respondents have not decided yet on staying at the university and $7.5\%$ have clear plan on leaving the academia definitely. 

Out of those who have strongly decided against second-cycle studies, over $45\%$ want to undertake additional courses; on the other hand, almost $40\%$ of them want to find employment as fast as possible.

Let us focus for a moment on the arguments which students gave to support their decision in case they wanted to continue their path outside of pursuing the mathematical Masters degree offered by their respective Alma Mater. Due to the latter being a multiple choice question, we first take a look at the number of times the selected reason was provided. Then we further analyse the relationship between those reasons and other answers they have provided in the questionnaire. Below we present the table with the responses from 82 students (after cleaning the data):

\begin{table}[H]
\begin{tabular}{ |p{0.75\textwidth}|p{0.125\textwidth}|p{0.125\textwidth}| } \hline
Statement & Count & Percentage\\\hline
I see no point in continuing studying mathematics from the perspective of my professional career. &$ 37 $& $ 45.12 \%$ \\\hline
Mathematical studies failed to meet my expectations regarding the level of scientific/knowledge standards. &$ 29 $& $ 35.37\%$ \\\hline
Mathematical studies are too time-consuming for me. &$ 24 $& $ 29.27\%$ \\\hline
The studies turned out to be too stressful for me. &$ 24 $& $ 29.27\%$ \\\hline
Mathematical studies failed to meet my expectations regarding the organization. &$ 22 $& $ 26.83\%$ \\\hline
I did not find mathematical studies interesting. &$ 19 $& $ 23.17\%$ \\\hline
I have no idea what to do after my studies. &$ 16 $& $ 19.51\%$ \\\hline
I want to change my university to more prestigious one. &$ 12 $& $ 14.63\%$ \\\hline
I want to change my university to less demanding one. &$ 7 $& $ 8.54\%$ \\\hline
My university does not offer specialization I would like to continue studying. &$ 7 $& $ 8.54\%$ \\\hline
I prefer to avoid certain professors from my university. &$ 4 $& $ 4.88\%$ \\\hline
Lack of part-time studies in mathematics. &$ 3 $& $ 3.66\%$ \\\hline
\end{tabular}
\caption{The reasons for not continuing the second-cycle studies on their original universities, provided by the students.}
\label{table:reasons}
\end{table}

The results can be disheartening from the standpoint of a lecturer. Most students do not see any real merit in further pursuing the education in the field of mathematics (this theory will be contrasted with the data from ELA database in subsequent part of this paper).
The second most commonly given reason is the one related to the level of education they receive. This is to some extent surprising when we compare it with the responses to questions about the readiness of their lecturers. Amongst other most commonly provided reasons were studies being either too stressful or too time-demanding. The latter should not be perceived as something unforeseen -- after all, the questionnaire has been conducted on people participating in full-time studies.

Utilizing Apriori algorithm we have discovered, that amongst most commonly associated reasons were:
\begin{itemize}
\item \textbf{A:} Mathematical studies failed to meet my expectations regarding the level of scientific/knowledge standards; \textbf{B:} Mathematical studies failed to meet my expectations regarding the organization -- this combination has over $60\%$ \textbf{A} $\to$ \textbf{B} confidence and over $80\%$ \textbf{B} $\to$ \textbf{A} confidence, while having nearly $22\%$ support.

\item \textbf{A:}  Mathematical studies failed to meet my expectations regarding the level of scientific/knowledge standards.; \textbf{B:} I see no point in continuing studying mathematics from the perspective of my professional career. -- this pair has also appeared together quite often (over $15\%$ of responses). It nears $45\%$ of \textbf{A} $\to$ \textbf{B} confidence and $37\%$ in reversed direction.

\item \textbf{A:} Mathematical studies failed to meet my expectations regarding the level of scientific/knowledge standards.; \textbf{B:} I did not find mathematical studies interesting.  -- This pairing occurs in $10\%$ of answers, with $31\%$ of confidence for \textbf{A}$\to$\textbf{B} and $47\%$ in other direction.

\item \textbf{A:} I did not find mathematical studies interesting.; \textbf{B:} Mathematical studies are too time-consuming for me. -- this matching has occurred in nearly $10\%$ of cases. It boasts $42\%$ \textbf{A} $\to$ \textbf{B} confidence and $36\%$ \textbf{B} $\to$ \textbf{A} confidence.

\item \textbf{A:} I see no point in continuing studying mathematics from the perspective of my professional career.; \textbf{B:} The studies turned out to be too stressful for me. -- such combination appears in $12\%$ of responses, with $28\%$ of \textbf{A} $\to$ \textbf{B} confidence and $45\%$ \textbf{B} $\to$ \textbf{A} confidence.
\end{itemize}

Additionally, amongst $34$ people who want to continue studying in different field, the most commonly provided reason for this change of heart was lack of perceived benefits from continuing the mathematical education ($68\%$ of responses). Over a third of them stated that they find the mathematical studies too stressful and almost $30\%$ of them found the topic of their studies wearisome.

\subsection{Comparison between retention group and leaving group}

By retention group we will consider students who would like to stay at their current university for second cycle of studies. Leaving group consists of remaining respondents. Let us start with the juxtaposition between the opinion of those two factions on general and specialized subjects. Due to lack of normality, we use Mann-Whitney U test for comparing the average similarity of their opinions.

\begin{table}[H]
\begin{tabular}{ |p{0.6\textwidth}|p{0.1\textwidth}|p{0.1\textwidth}|p{0.1\textwidth}|p{0.1\textwidth}|} \hline
Question & Mean~of retention group & Mean~of leaving group & p-value \\\hline
\textbf{General courses contain too much theory and too little practice.} & 2.89 & 3.36 & 0.0183\\\hline
General courses contain too much practice and too little theory. & 1.95 & 1.87 & 0.2826\\\hline
\textbf{Lecturers in general courses are well-prepared for classes.} & 4.1 & 3.74 & 0.0102\\\hline
\textbf{General courses prepared me for the expectations of the job market.} & 2.71 & 2.26 & 0.0009\\\hline
\textbf{General courses provided me with soft skills such as logical thinking and evaluating the consistency of reasoning, which will undoubtedly be useful in the job market.} & 3.94 & 3.56 & 0.0126\\\hline
There are too many computer science subjects within the general courses. & 2.3 & 2.4 & 0.7948\\\hline
There are too few computer science subjects within the general courses. & 2.61 & 2.63 & 0.9133\\\hline
\textbf{Studies should be integrated programs (5-year programs leading to a Master's degree)}. & 2.72 & 2.23 & 0.0212\\\hline
First-cycle mathematics studies should end with an engineering degree, even if
 this means the necessity of studying technical subjects (e.g. technical drawing, physics). & 2.78 & 2.81 & 0.8263\\\hline
\end{tabular}
\caption{Comparison of students' opinions on general courses. Bold font was used to point out statistically significant differences.}\label{table:gentable}
\end{table}

A very succinct summary of the above Table \ref{table:gentable} would be -- students who are reluctant towards staying at their current university for Masters degree in mathematics consider the general courses as much less useful for their future careers in comparison with retention group. They consider general courses to be too theoretical and are quite strongly against the idea of 5-year Master's studies with no Bachelor title after three years. Similar conclusions can be drawn from the differences in opinions regarding the specialized courses.

In the following Table \ref{table:spectable} we present the analogous comparison for the opinions on specialized courses. While the overall ordering of the two groups is preserved -- the leaving group remains more critical of lecturer preparation, perceived job-market relevance and soft-skills acquisition -- two shifts are worth emphasising. First, the theory-vs-practice disagreement that was sharp for general courses evaporates for specialized ones (p-value $0.5252$), which suggests that the over-abstraction complaint is directed at the general curriculum rather than at academic mathematics as such.
 
Second, the single largest between-group gap in the whole comparison concerns satisfaction with the chosen specialization itself (p-value $0.0003$), accompanied by a significantly stronger demand from the leaving group for a larger share of specialized courses in the curriculum (p-value $0.0039$). Taken together, these two observations hint that the leaving group's dissatisfaction is less about the nature of the specialized content and more about insufficient access to it and limited fit with the specialization menu on offer -- a point we return to in the recommendations section.

Another important observation from both tables is that opinion on computer-science subjects remains quite consistent in relationship to both specialized and general parts of the curriculum. 

\begin{table}[H]
\begin{tabular}{ |p{0.6\textwidth}|p{0.1\textwidth}|p{0.1\textwidth}|p{0.1\textwidth}|p{0.1\textwidth}|} \hline
Question & Mean~of retention group & Mean~of leaving group & p-value \\\hline
Specialized courses contain too much theory and too little practice. & 2.71 & 2.84 & 0.5252\\\hline
Specialized courses contain too much practice and too little theory. & 2.28 & 2.27 & 0.7122\\\hline
\textbf{Lecturers in specialized courses are well-prepared for classes.} & 4.14 & 3.66 & 0.0016\\\hline
\textbf{Specialized courses prepared me for the expectations of the job market.} & 3.3 & 2.88 & 0.0146\\\hline
\textbf{Specialized courses provided me with soft skills such as logical thinking and evaluating the consistency of reasoning, which will undoubtedly be useful in the job market.} & 3.84 & 3.42 & 0.0242\\\hline
There are too many computer science subjects within the specialized courses. & 2.28 & 2.31 & 0.6504\\\hline
There are too few computer science subjects within the specialized courses. & 2.65 & 2.74 & 0.6813\\\hline
Specialized courses prepared me to some extent for the job market, but I need to catch up on many things on my own. & 3.41 & 3.61 & 0.0827\\\hline
\textbf{There are too few specialized courses compared to general courses.} & 3.12 & 3.61 & 0.0039\\\hline
There are too many specialized courses compared to general courses. & 2.12 & 2.01 & 0.3034\\\hline
The selection of thesis topics related to my specialization is suitable for students' needs. & 3.14 & 2.84 & 0.0797\\\hline
\textbf{I am satisfied with the specialization I chose.} & 4.19 & 3.58 & 0.0003\\\hline
\end{tabular}
\caption{Comparison of students' opinion on specialized subjects. Bold font was used to point out statistically significant differences.}\label{table:spectable}
\end{table}

Lastly, the thesis-topic suitability nearly reaches significance (with observed p-value $0.0797$), pointing softly in the same direction as our previous observations: the leaving group feels the elective/specialization offer does not fully accommodate their interests.

\subsubsection{A Benjamini–Hochberg False-Discovery Rate correction}

When testing many hypotheses at once, simply using a p-value limit value $0.05$ can often lead to false positives. The Benjamini–Hochberg (BH) procedure \cite{Benjamini1995Controlling,Barber2017The} is a simple rule to control the false discovery rate (FDR) in such scenarios. For the discussed multiple comparisons across the 21 between-group tests reported in this section, we applied the BH procedure at $\alpha = 0.05$. The results are summarized in the subsequent table -- all of them survive the subsequent verification, even if barely.
Exact results are presented in the Table \ref{table:postbhcorec} below (which includes only statistically significant findings from the previous section).

\begin{table}[H]\hspace{-2cm}
\begin{tabular}{|p{0.5\textwidth}|p{0.1\textwidth}|p{0.1\textwidth}|p{0.1\textwidth}|p{0.1\textwidth}|p{0.1\textwidth}|}
\hline
Question & Retention group mean & Leaving group mean  & Raw p-value & FDR-corrected p-value & Rank-biserial correlation \\\hline
I am satisfied with the specialization I chose. & 4.19 & 3.58 & 0.0003 & 0.0054 & -0.31\\\hline
General courses prepared me for the expectations of the job market. & 2.71 & 2.26 & 0.0009 & 0.0099 & -0.275\\\hline
Lecturers in specialized courses are well-prepared for classes. & 4.14 & 3.66 & 0.0016 & 0.0114 & -0.265\\\hline
There are too few specialized courses compared to general courses. & 3.12 & 3.61 & 0.0039 & 0.0204 & 0.249\\\hline
Lecturers in general courses are well-prepared for classes. & 4.1 & 3.74 & 0.0102 & 0.0427 & -0.208\\\hline
General courses provided me with soft skills such as logical thinking and evaluating
 the consistency of reasoning, which will undoubtedly be useful in the job market. & 3.94 & 3.56 & 0.0126 & 0.0439 & -0.212\\\hline
Specialized courses prepared me for the expectations of the job market. & 3.3 & 2.88 & 0.0146 & 0.0439 & -0.2\\\hline
General courses contain too much theory and too little practice. & 2.89 & 3.36 & 0.0183 & 0.048 & 0.203\\\hline
Studies should be integrated programs (5-year programs leading to a Master's degree). & 2.73 & 2.23 & 0.0212 & 0.0494 & -0.199\\\hline
\end{tabular}
\caption{A re-examination of statistical significance of  students' opinion on specialized subjects.}\label{table:postbhcorec}
\end{table}

In subsequent sections the FDR corrections will be included along the performed analyses. 

For measuring effect sizes we have decided to use rank-biserial correlation \cite{Tapio2025The, Willson1976Critical}, which is a convenient way to quantify how strongly a dichotomous group variable is associated with an outcome of ordinal/ranked nature. Values close to 0 indicate little separation of the groups' ranks; larger absolute values indicate stronger separation. Benchmarks are not universally agreed on -- certain papers refer to effect as negligible if the coefficients' absolute value does not exceed 0.10-0.15. Moderate effects tend to be rated between 0.3 and 0.6 in terms of the absolute value of this coefficient and when the effect size exceeds 0.7, most sources agree that the groups are well-separated on the outcome.\footnote{Examples of  similar interpretations for effect size or association strength can be found in education field as well, see, e.g. \cite{Anabo2023CORRELATES,Pangesti2025Application}.} Since such thresholds are often domain-specific and somewhat arbitrary, we limit ourselves to reporting the coefficient values directly rather than assigning qualitative labels.

\subsection{Demographics and background impact}\label{subsec:genderr}

We begin with comparing how the opinions on general and specialized courses differ between genders -- a factor that remains one of the important variables in retention in mathematically-demanding STEM fields, especially given the influence of social stereotypes (see recent papers of Almukhambetova et al. \cite{Almukhambetova2023LeakyPipeline,Almukhambetova2023RetentionSTEM} and references therein, as well as \cite{hill2010so,silbey2016so}). In Table \ref{table:gender} we report only those questions, where p-value was below $0.2$ for the sake of brevity.\footnote{Due to the small sample size of other responses, Table \ref{table:gender} includes only male and female respondents.} Only one finding survives FDR correction, with second one remaining almost statistically significant. It seems that female students not only have stronger opinion on the topic of engineering degree, but they also seem more aware on the work necessary to do by themselves in order to be more recognized on the job market.

\begin{table}[H]
\hspace{-2cm}
\begin{tabular}{|p{0.5\textwidth}|p{0.1\textwidth}|p{0.1\textwidth}|p{0.1\textwidth}|p{0.1\textwidth}|p{0.1\textwidth}|}
\hline
Question & Female mean & Male mean  & Raw p-value & FDR-corrected p-value & Rank-biserial correlation \\
\hline
\textbf{First-cycle mathematics studies should end with an engineering degree, even if
 this means the necessity of studying technical subjects (e.g. technical drawing, physics).} & \textbf{3.49} & \textbf{2.22}  & \textbf{0.0001} & \textbf{0.0020} & \textbf{-0.4718} \\\hline
Specialized courses prepared me to some extent for the job market, but I need to catch up on many things on my own. & 3.91 & 3.33 & 0.0056 & 0.0590 & -0.3261 \\\hline
General courses contain too much theory and too little practice. & 3.72 & 3.09 & 0.0482 & 0.3374 & -0.2362 \\\hline
Specialized courses provided me with soft skills such as logical thinking and evaluating the consistency of reasoning, which will undoubtedly be useful in the job market. & 3.23 & 3.53 & 0.1089 & 0.4743 & 0.1897 \\\hline
General courses contain too much practice and too little theory. & 1.61 & 2.00 & 0.1129 & 0.4743 & 0.1809 \\\hline
General courses prepared me for the expectations of the job market. & 2.09 & 2.38 & 0.1552 & 0.5431 & 0.1674 \\\hline
Studies should be integrated programs (5-year programs leading to a Master's degree). & 2.35 & 2.00 & 0.1952 & 0.5665 & -0.1530 \\\hline
\end{tabular}
\caption{Comparison of opinions on learning aspects at their universities between genders.}\label{table:gender}
\end{table}

As far as the financial self-positioning is concerned, it does not seem to impact the decision to change their plans regarding second cycle of studies.

\begin{figure}[H]
\centering
\includegraphics[scale=0.55]{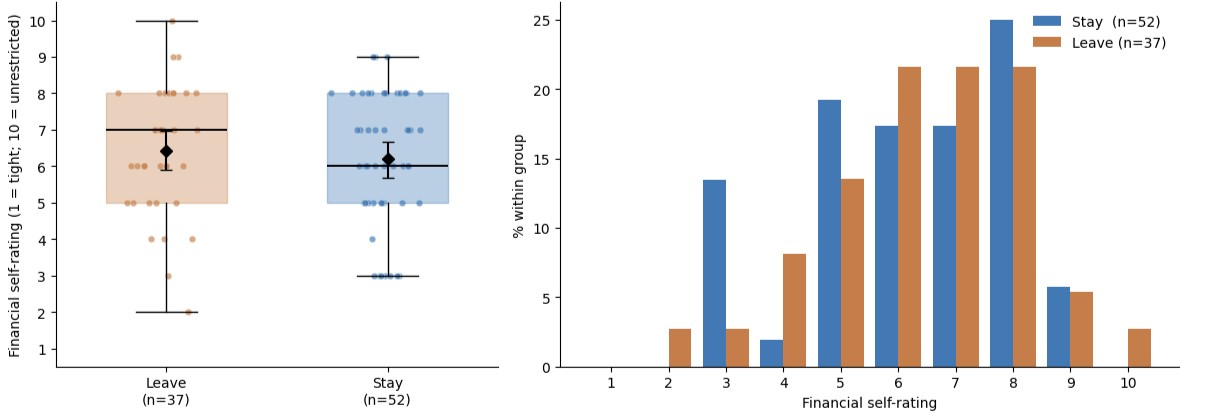}
\caption{The comparison of self-assessed financial situation for each group.}
\end{figure}

The mean rating of staying group was $6.21$, while the leaving students on average rated their financial standing on $6.43$. Mann-Whitney U test was applied due to lack of normality with two-sided alternative hypothesis, which yielded p-value $0.651$. Very small differences in distribution were observed (staying group sd: $1.786$ vs leaving group: $1.741$), but they are not statistically significant according to the Levene test (p-value $0.661$). The fact that the financial situation is evaluated through self-perception rather than by crisp, monetary lens make these results harder to compare with findings of other authors, who often claim financial strain as one of the most important factors which contribute towards the decision of leaving the university (see \cite[Figure 2]{nieuwoudt2021whystudents}, where financial hardship contributed to over $44\%$ of the MCQ answers as well as \cite{heublein2007ursachen, Geisler2018}, where this was the 4th most common main reason for leaving).

\begin{figure}[H] 
 \centering
 \includegraphics[height=0.6\textheight]{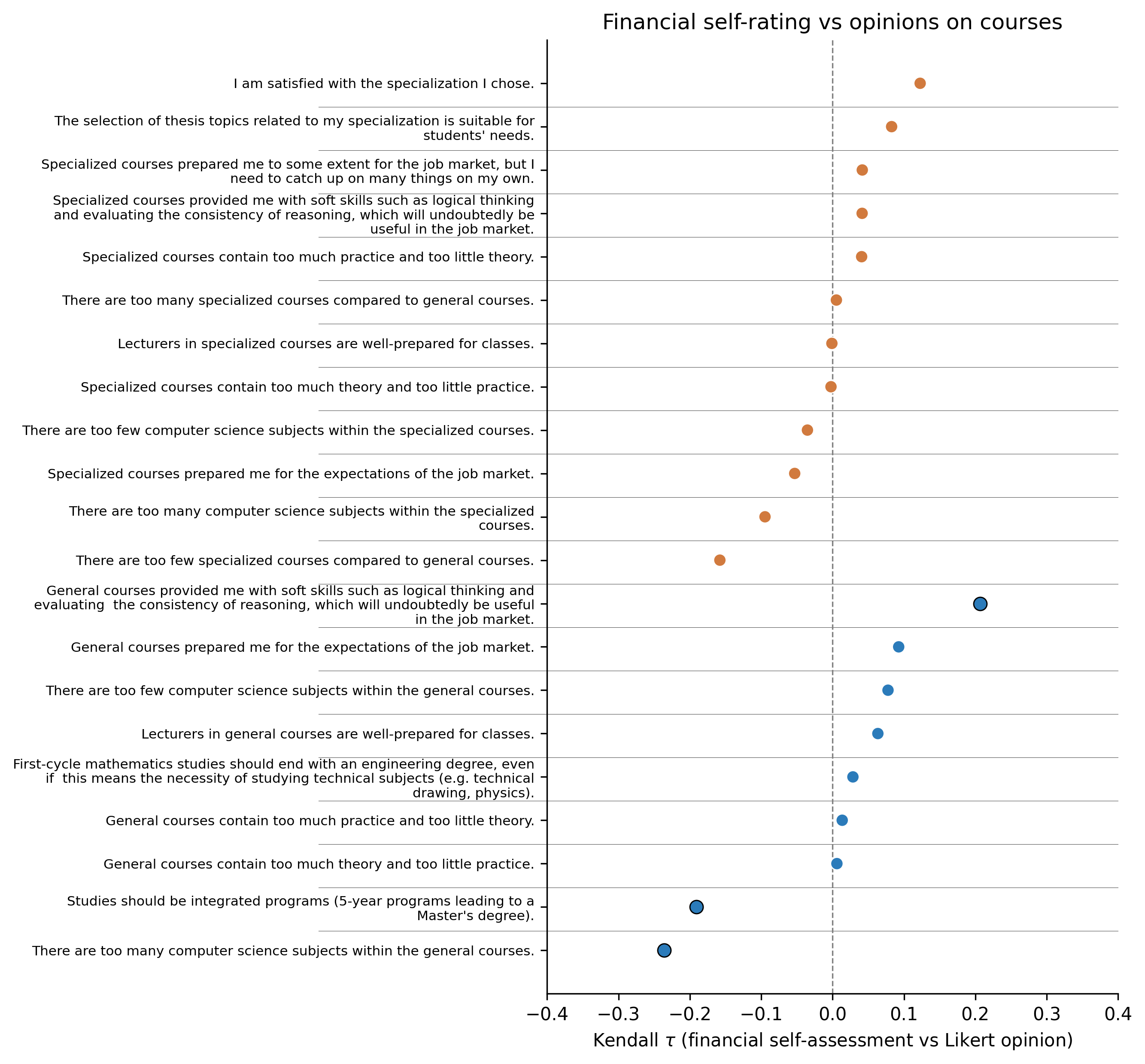}
 \caption{Correlation analysis between financial self-assessment and investigated opinions on specialized (orange) and general (blue) courses. The black outline around the point indicates statistical significance prior to FDR correction.}
\label{fig:corrtree}
\end{figure}

Students from lower-income socioeconomical backgrounds are more likely to perform a strict cost-benefit rationale calculus \cite{Mller2022}. We try to engage this idea by looking at the correlation between financial situation and the opinion on the courses offered by the university. For sake of brevity in Table \ref{table:corr} we have included only statistically significant correlations (remaining ones were both low and insignificant with absolute magnitudes below $0.15$). Full picture of the correlations grouped according to the question categories can be seen in Figure \ref{fig:corrtree}. 
It is important to point out that no correlation remains significant after FDR correction.

\begin{table}[H]
\begin{tabular}{|p{0.5\textwidth}|p{0.1\textwidth}|p{0.1\textwidth}|p{0.1\textwidth}|p{0.1\textwidth}|}
\hline
Question & Mean opinion & Kendall $\tau$ & Raw p-value & FDR-adjusted p-value \\
\hline
There are too many computer science subjects within the general courses. & 2.3708 & -0.2359 & 0.0058 & 0.1224 \\\hline
General courses provided me with soft skills such as logical thinking and evaluating
 the consistency of reasoning, which will undoubtedly be useful in the job market. & 3.5843 & 0.2066 & 0.0166 & 0.1742 \\\hline
Studies should be integrated programs (5-year programs leading to a Master's degree). & 2.2472 & -0.1910 & 0.0260 & 0.1817 \\\hline
 \end{tabular}
 \caption{Correlations between found opinions on courses and financial self-assessment.} \label{table:corr}
 \end{table}
 
Concluding this section, we also notice that working alongside studies is also commonly indicated as one of the drop-out predictors (again, see \cite{nieuwoudt2021whystudents}). In our research, the students who work alongside studies are, surprisingly, slightly more willing to stay on average, but the difference is not statistically significant (prop test statistic $z = 0.827$,  p-value: $0.4081$). The two sided Mann-Whitney $U$ test did not confirm the significance of difference between average working hours between the working retention subgroup ($23.38$h per week on average)  and the working students who do not intend to continue studying on their original university ($15.50$h) -- p-value $0.1435$, rank biserial coefficient $-0.394$.

 \begin{figure}[H] 
 \centering
 \includegraphics[width=0.6\textwidth]{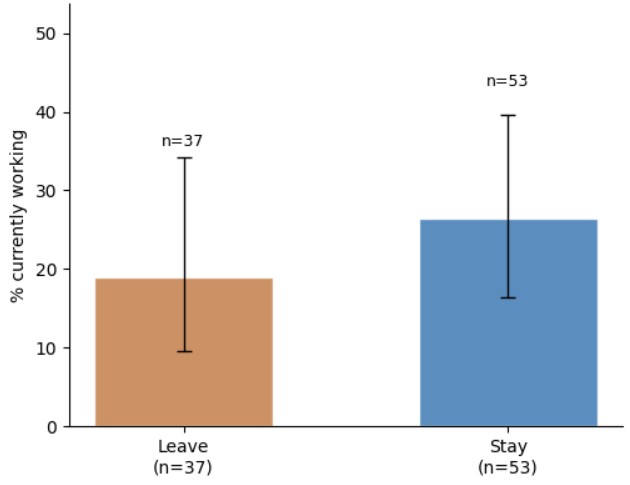}
 \caption{Comparison of working status between groups.}\label{fig:wgroup}
 \end{figure}

This can be motivated partially by the fact that people who have started working see the value of the education they receive. This statement will  be discussed in the Section \ref{sec:elaa}, where we perform a brief analysis of graduates information collected via Polish Graduate Tracking Systemm.

\subsection{Academic performance impact on staying}

Various sources indicate that the academic performance has significant impact on retention, especially in the early stages of studies \cite{Kocsis2024Factors}. We are more inclined towards the conclusions postulated by Espinoza et al. \cite{espinoza2023whydo} who argue that the relationship between GPA and withdrawal from studies is more complex than it initially meets the eye. It is a rather widespread myth that most people who leave the university simply performed too poorly to stay, as there are observable deviations from this rule (which are also noticeable in this study). In particular, factors related to purely demographic factors, academic performance and achievements contributed in total to just 12\% of the total variance of degree completion (see \cite{nieuwoudt2021whystudents, IJESL2237} and references therein). On the other hand, the notion of mathematical self-concept (explained quite intuitively in \cite{Geisler2023,rach2019self}, also see references pointed therein) has a positive relationship with persistence and other achievement related features which, in turn, have positive impact on students' retention. The GPA is obviously not directly translatable to the mathematical self-concept, although this is the closest fit we have available at our hands. 

In general, our findings in this topic are best summarized by the following picture:

\begin{figure}[H]
\centering
\includegraphics[scale=0.55]{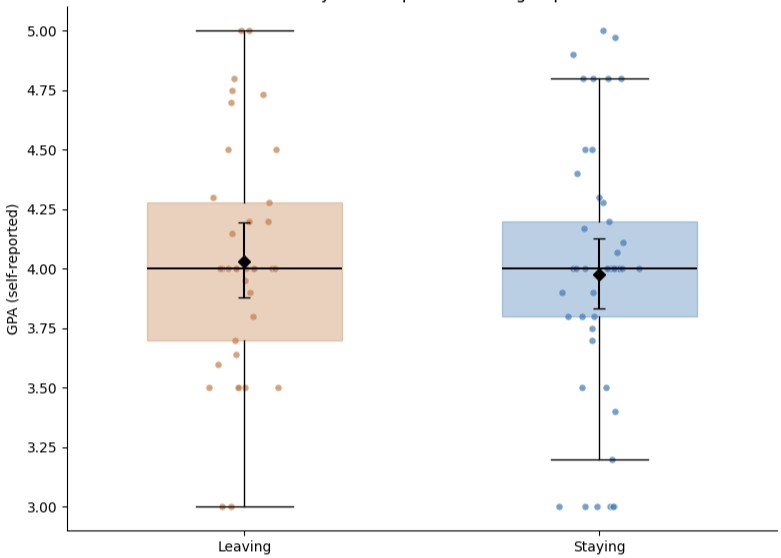}
\caption{The comparison of self-reported GPA from the last year of studies for each group.}
\end{figure}

The group which has decided on continuing studies had slightly lower self-reported GPA ($3.977$ vs $4.032$), but the difference is nowhere near being statistically significant (p-value $0.797$ for Mann-Whitney U test). The variances in GPAs also do not differ significantly (standard deviations for remaining group: $0.540$ vs    leaving group: $0.490$, Levene test p-value $0.777$). 

Concluding, this study does not indicate that GPA is directly related to the phenomenon of abandoning the alma mater after the first cycle of studies, albeit additional, independent confirmation of it would be appreciated.

%
%

\section{Comparing the graduates -- a brief analysis}\label{sec:elaa}

While the questionnaire results provide insight into students' declared motivations and intentions, they do not directly answer whether abandoning second-cycle studies is economically justified from the graduates' perspective. To partially address this issue, we complement the survey-based analysis with administrative labour-market data obtained from the ELA system. 

The Polish nationwide graduate tracking system \emph{Ekonomiczne Losy Absolwent\'ow} (ELA, Economic Fates of Graduates), publicly available at \url{www.ela.nauka.gov.pl}. ELA is one of the most innovative graduate monitoring systems in Europe and, unlike survey-based approaches used in many countries, relies entirely on administrative records. Two registers feed the system: the Polish Social Insurance Institution (\emph{Zak{\l}ad Ubezpiecze\'n Spo{\l}ecznych}, ZUS), which provides monthly data on employment status, contract type, and contributions paid on behalf of insured individuals, and POL-on, the national register of higher education, which supplies information on the graduate's field of study, institution, degree level, and year of graduation. Place of residence at the county (\emph{powiat}) level is identified via the TERYT register, which makes it possible to contextualize each graduate's labour market outcomes against statistics on the local labour market published by Statistics Poland (GUS).

Because data are drawn directly from administrative registers rather than collected through respondent contact, the system covers the entire population of higher education graduates (with the exception of certain military academy graduates), eliminating sampling error and the recall, rounding, and social desirability biases typical of self-reported data.
Graduates are distinguished only by the assigned identifier; no personal data are stored, so anonymity is fully preserved. To further protect privacy, indicators are not reported for cohorts smaller than 10 graduates, and breakdowns involving sub-groups of fewer than 3 individuals are suppressed. Several limitations of administrative data should nonetheless be borne in mind. ZUS records do not capture earnings from civil-law contracts (\emph{umowa o dzie{\l}o}, \emph{umowa zlecenie}) signed with students under the age of 26, work performed abroad, informal employment,\footnote{For example, private mathematical tuitions which are quite common way for a STEM-oriented students to gain some extra money during their studies (in most cases, without reporting it to ZUS).} or activity covered by the separate farmers' insurance scheme (KRUS); income from self-employment is also unobserved, as self-employed individuals typically declare the minimum permitted contribution base regardless of actual earnings. POL-on, in turn, does not record academic grades or final diploma marks. The tenth edition of ELA, used here, covers graduates from the years 2014--2023 and follows the 2014--2018 cohorts over a full five-year post-graduation window.

Before proceeding with the subsequent analyses, we note that we intentionally avoided labelling curves corresponding to individual universities. Although an interested reader could, with some effort, trace these results in the ELA database, our aim was to prevent the paper from being interpreted as an evaluation of specific institutions.

\subsection{Earning premium of Masters over Bachelors}

In Subsection \ref{subsec:retentionbetween} we have indicated that we were going to test the \textit{lack of benefits} theory regarding continuation of the second-cycle mathematical studies, which was the most commonly pointed out reason in Table \ref{table:reasons}. Let us begin with the following spaghetti plot, based on the data collected by ELA. The earnings reported here are combined monthly earnings from all sources in year following the graduation.

\begin{figure}[H]
\centering
\includegraphics[scale=0.55]{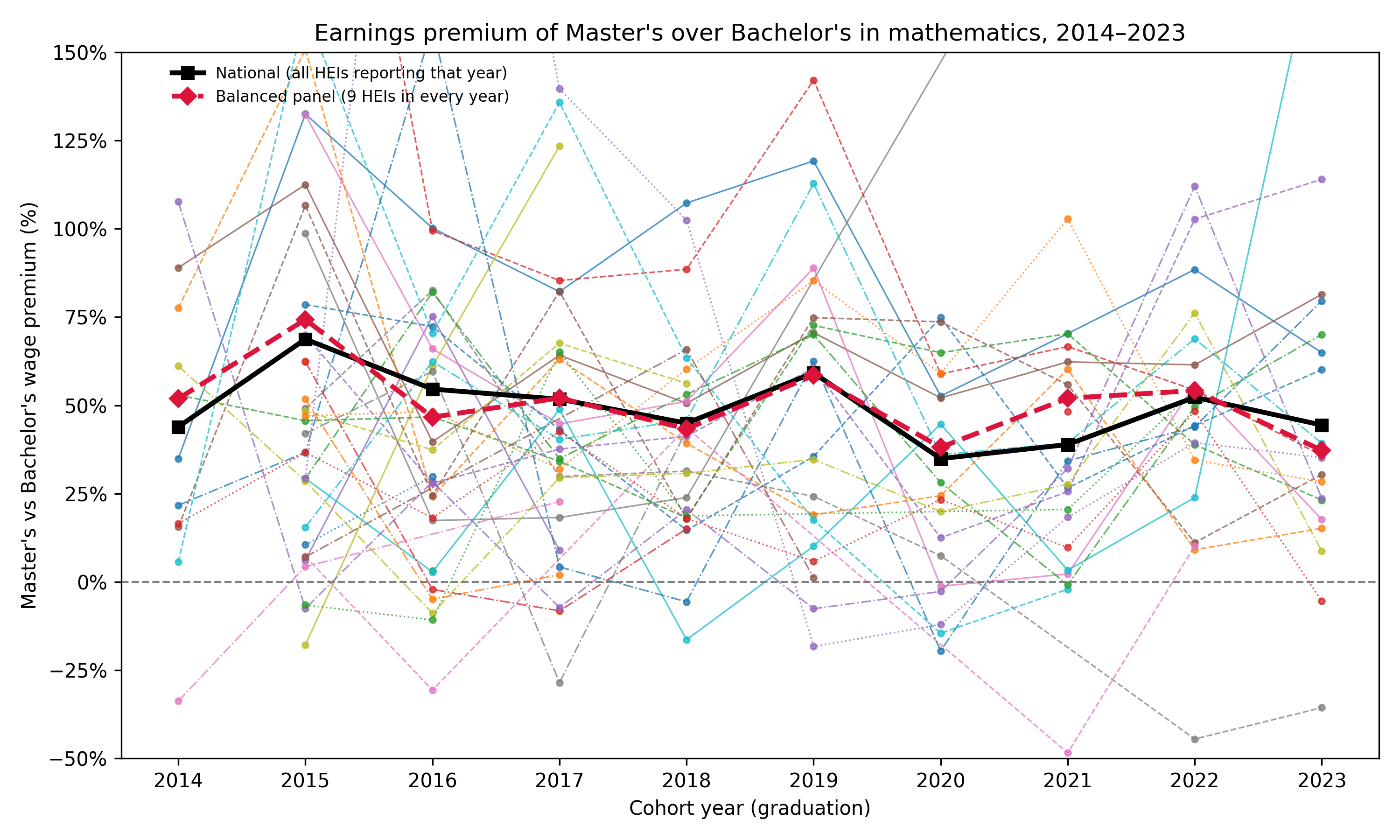}
\caption{The earnings advantage provided by the Masters degree in the first year following the graduation.}\label{fig:premiumhei}
\end{figure}

Several per-HEI lines are interrupted in the figure (as well as in the Figures to come: this reflects the fact that ELA suppresses indicator cells whenever the underlying cohort of mathematics graduates falls below the institute's disclosure threshold, which is common for smaller programmes and certain cycles in certain years. The pattern is therefore an artefact of reporting rules and not a definitive evidence of programme closure.

The set of reporting HEIs shifts from year to year, the headline national line could, in principle, move simply due to compositional changes. The dashed crimson line restricts the weighted mean to the subset of HEIs that report in every year of the period, and tracks the unrestricted national line closely — confirming that the temporal pattern is genuine rather than driven by entry and exit of institutions in the panel.
Despite certain outliers, the general tendency remains rather stable, granting a stable bonus earnings to the Masters, oscillating around $45$-$50\%$ for all years reported in the database.

One could make a (rather valid) point that such comparison is biased -- after all, people who abandon studies gain experience in the two years following the graduation, whereas their academic colleagues fall back in that aspect. In subsequent plot we compare the earnings (from all accounted sources) between the following groups.

\begin{itemize}
\item Masters in the year following their graduation;
\item Bachelors in the third year following their graduation;
\end{itemize}

The subsequent pair of plots prove that the intuition about the first-cycle studies combined with acquired experience being worth more than additional degree seems correct. What needs to be pointed out, however, is the short-term nature of this comparison, limited by the fact that the ELA system offers data which covers up to 5 years of the graduates' careers. Shorter range of cohort years analysed therein is due to the data limitations.

\begin{figure}[H]
\centering
\includegraphics[scale=0.55]{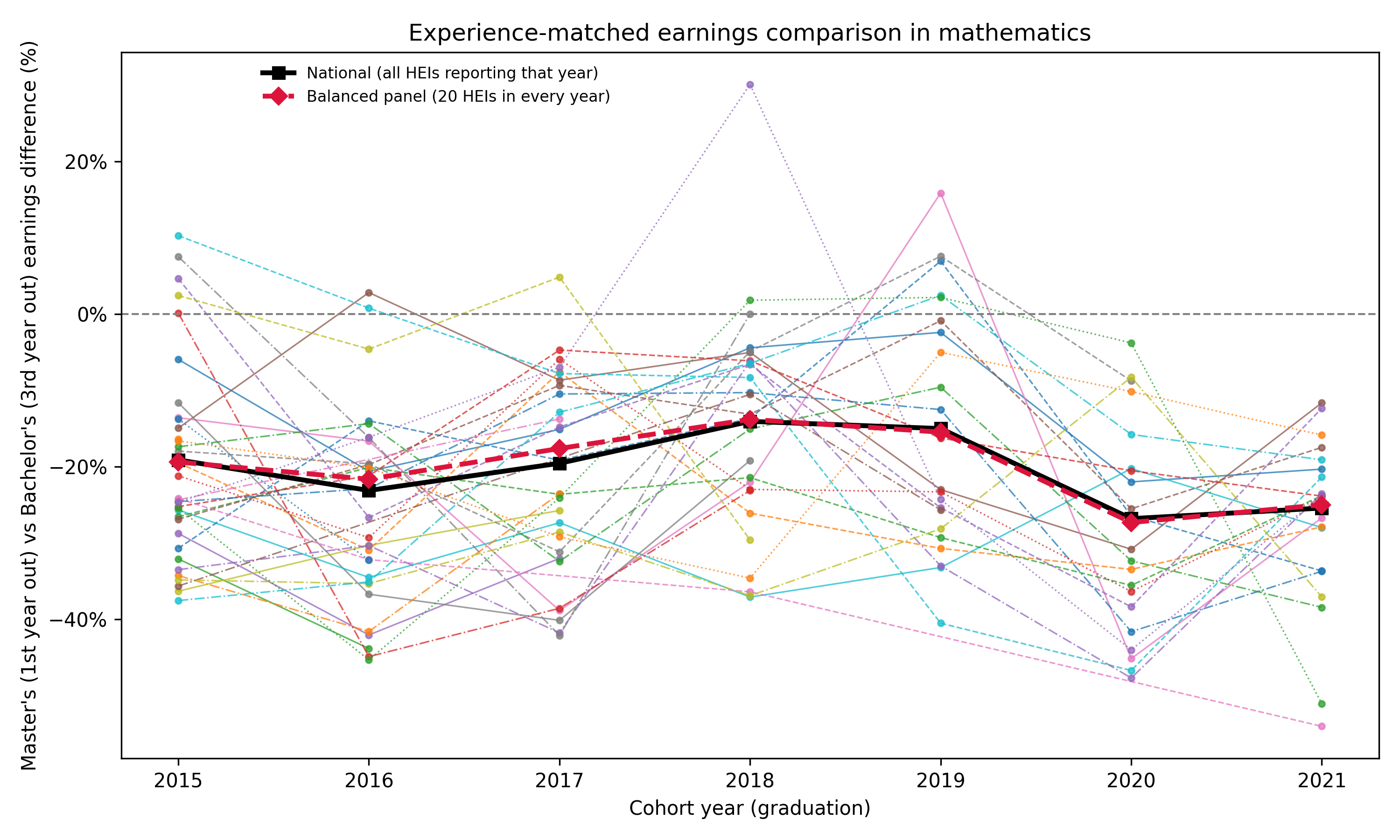}
\caption{The earnings ratio between students staying for Masters degree and those who spent this time at work -- Master's at year 1 vs Bachelor's at year 3 post-graduation.}\label{fig:premiumhei3}
\end{figure}

\begin{figure}[H]
\centering
\includegraphics[scale=0.55]{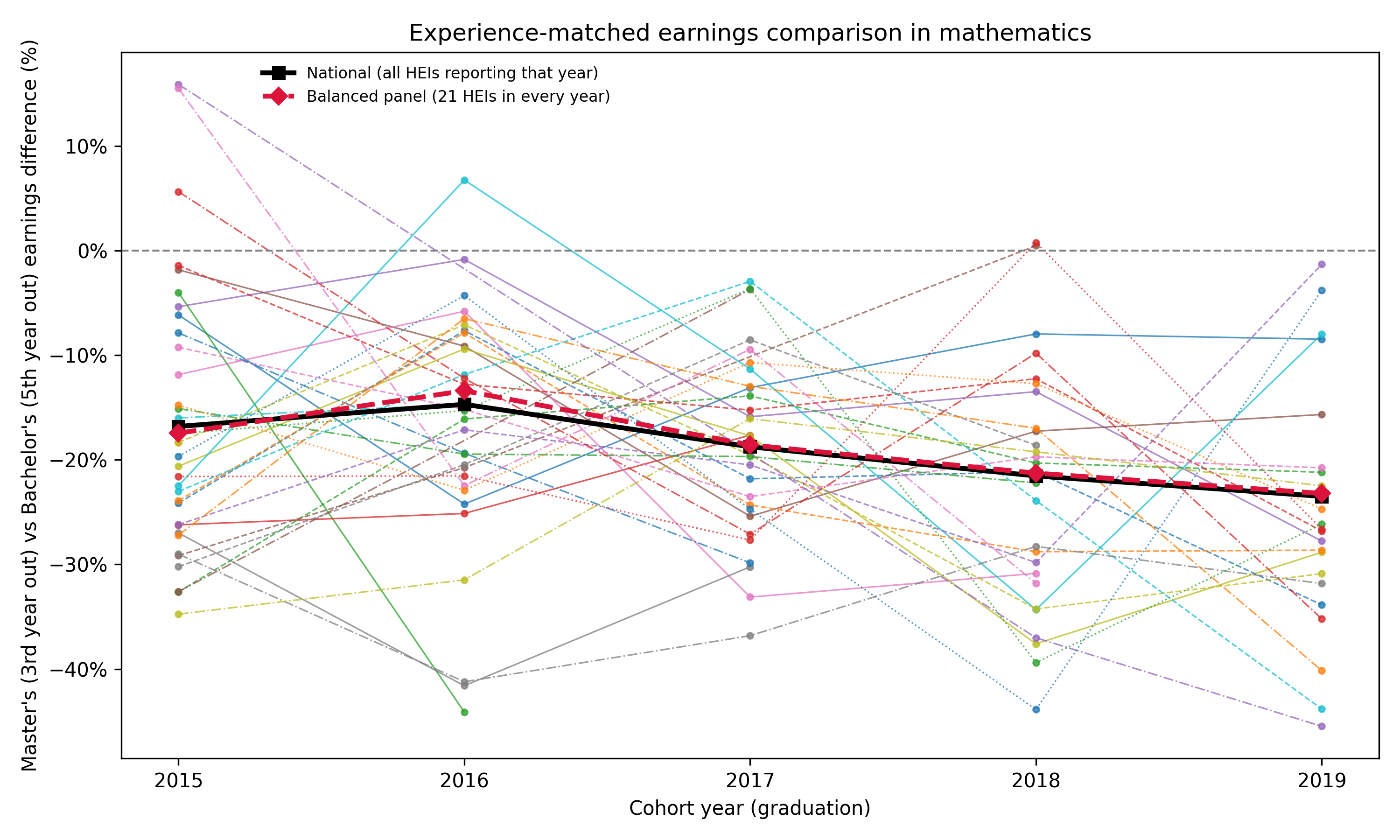}
\caption{The earnings ratio between students staying for Masters degree and those who spent this time at work -- Master's at year 3 vs Bachelor's at year 5 post-graduation.}\label{fig:premiumhei5}
\end{figure}

This trend might seem ridiculous from the educators' perspective, but completely reasonable when perceived through the job market lenses. After all, a person who has obtained Bachelor in mathematics has to operate mathematical formalism quite well. Instead of deepening that understanding and building subsequent layers of mathematical knowledge, the practical job experience obtained in this time can sharpen and shape their skills in direction desired by the employers. As sad as it might be, no amount of time in academia can completely replace the expertise granted by actual experience -- but that phenomenon can also be viewed as a call for the action for the people responsible for creating syllabuses and general organization of the second-cycle studies in mathematics. Coordination of the learning process with the job market should be at least thoroughly considered by the decision-makers in this field.

\subsection{Time to first employment comparison}

On average, people who managed to complete their second-cycle studies have much easier entry to the job-market, as displayed in the following picture:

\begin{figure}[H]
\centering
\includegraphics[scale=0.55]{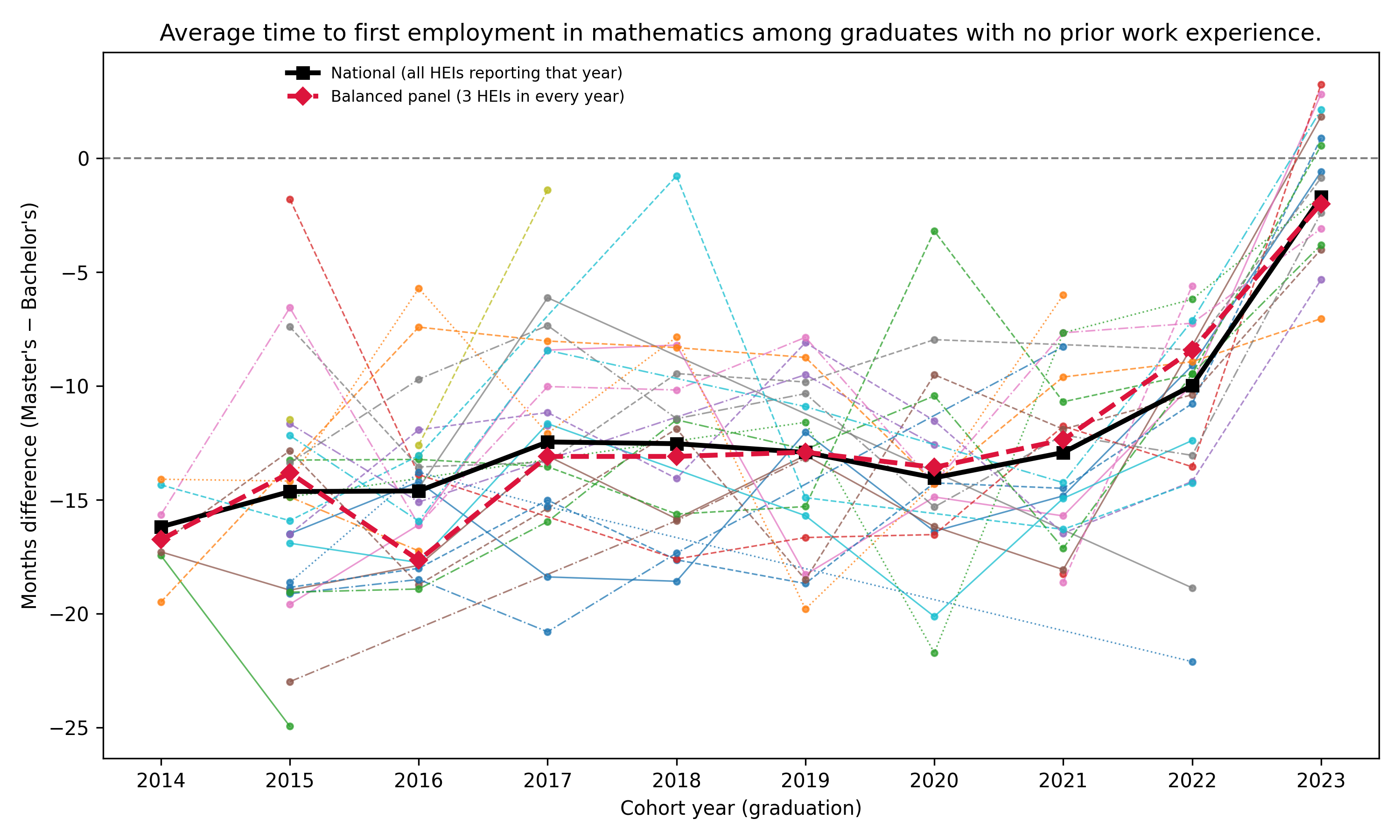}
\caption{The difference between time spent on finding employment after graduation (without prior job experience) for Masters and Bachelors.}\label{fig:months}
\end{figure}

Negative values indicate better performance of respective year Masters, while positive suggest that Bachelors do better on average. As we can see, the Masters degree used to grant strong advantage in this aspect. However, from the 2020
cohort onwards this advantage visibly contracts. While a single explanation is unlikely to capture the entirety of this shift, a plausible contributing factor is the proliferation of expedited recruitment channels that gained prominence during the pandemic period --- in particular the streamlined, often algorithmically-mediated application tracks promoted by
large job-listing platforms. Such services lower the marginal cost of applying for less experienced candidates and thereby compress the search-to-hire interval at the Bachelor's end of the labour market disproportionately. This drift can also be partially shaped by the ever-changing job market expectations.

\subsection{Longitudinal context: realized retention in mathematics from ELA dataset}

While our questionnaire captures the second-cycle intentions of a single cohort of students surveyed in 2024, the OPI-PIB ELA system makes it possible to
contextualize those intentions against nearly a decade of \emph{realized} retention behaviour.

\begin{figure}[H]
\centering
\includegraphics[width=\textwidth]{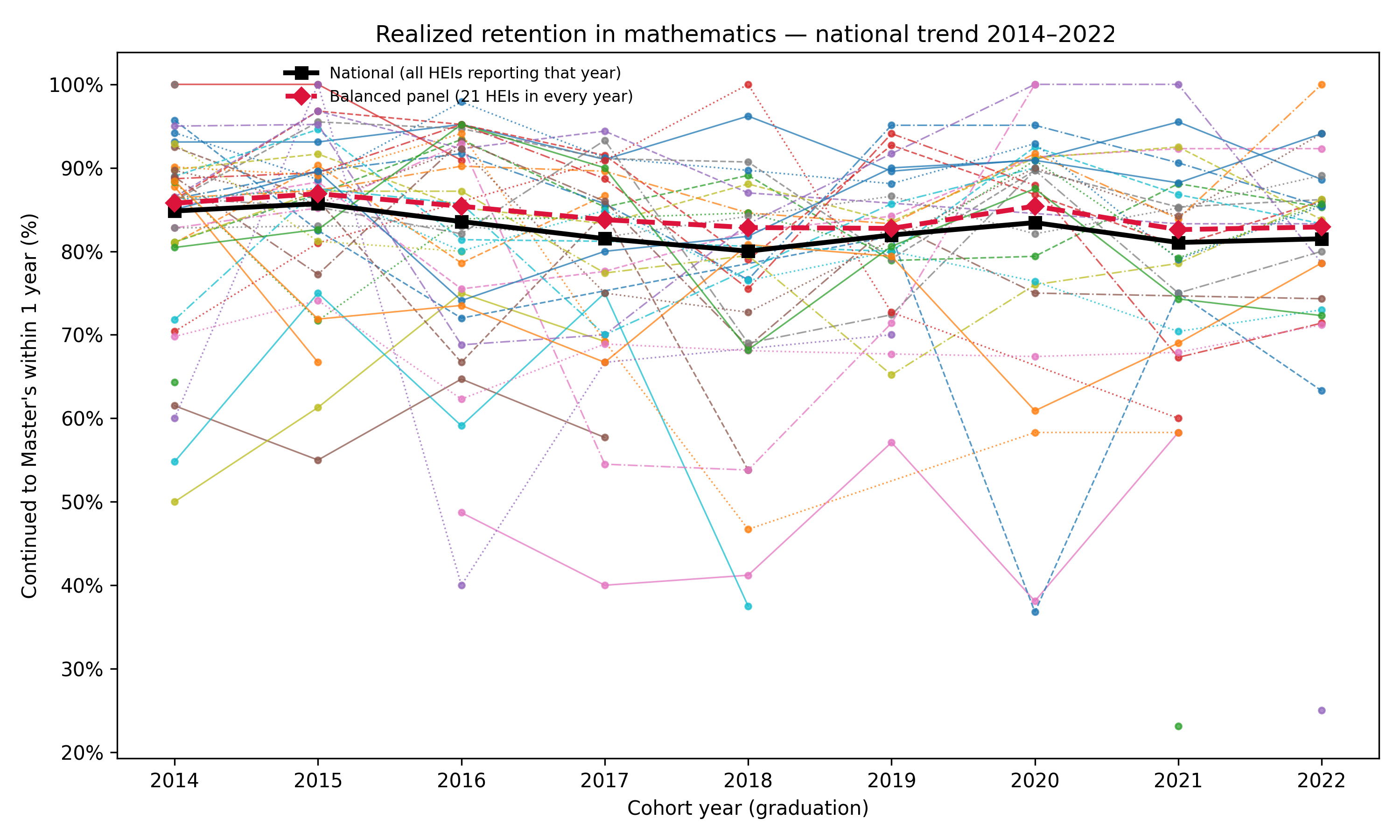}
\caption{The retention behaviour for the mathematical studies (allowing possible different university).}\label{fig:retenhei}
\end{figure}

\begin{figure}[H]
\centering
\includegraphics[width=\textwidth]{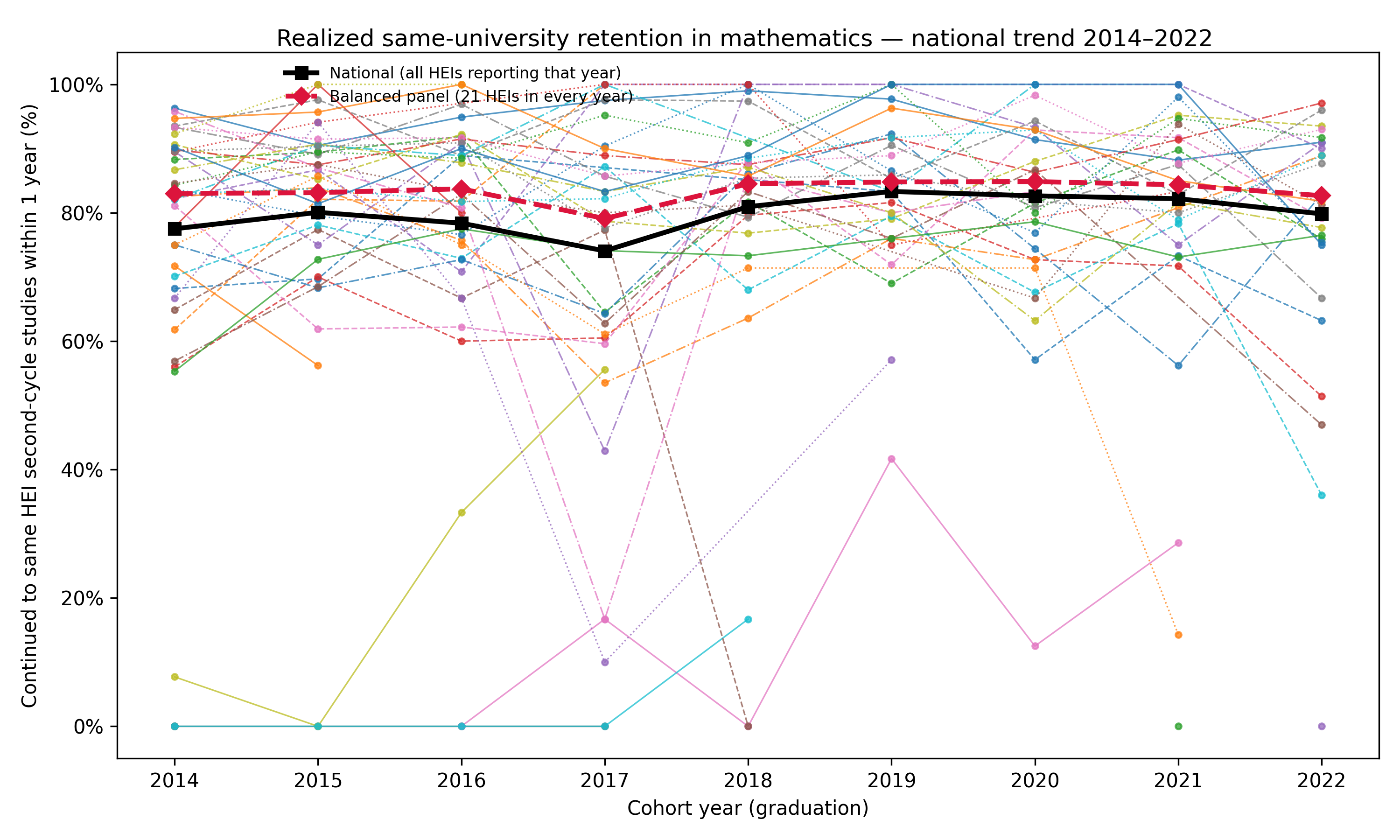}
\caption{The retention behaviour for the mathematical studies (within same university).}\label{fig:retenhei2}
\end{figure}

Figures~\ref{fig:retenhei} and \ref{fig:retenhei2} report, for every graduation cohort between 2014 and 2022, the percentage of mathematics
first-cycle graduates who enrolled in second-cycle studies within one year of obtaining their Bachelor's degree in mathematics, applied mathematics and their English-language counterparts.

Once again, the black curve represents the national rate, weighted by the number of graduates in each (HEI, cohort) cell; the dashed crimson curve recomputes the same weighted mean restricted to the balanced
panel of institutions that report data for every year in the period, serving as
a robustness check against compositional drift. The two aggregate curves track
each other closely, indicating that the headline pattern is not driven by entry
or exit of HEIs from the reporting set.

Both Figures present stable tendency at the national level, accompanied by considerable and chaotic cross-institutional
dispersion: the difference between the highest- and lowest-retention HEI in any
given year regularly exceeds $55$ percentage points. This longitudinal evidence
strengthens two observations from our 2024 cross-sectional survey. First, the stay-rate stated by our respondents (both regarding the same-university and same-subject retention) is not entirely consistent with
the realized rates documented by ELA in nearby cohorts. Nevertheless, the ratio of students who choose to stay for mathematics on the same university to those who wish to study mathematics in general remains relatively similar. Second, the persistence of large between-HEI gaps over almost a decade indicates that retention is, at least to some extent, an institution-level phenomenon -- a finding that any system-wide intervention motivated by our recommendations in Section~\ref{sec:recommend} should take into account.

\section{Study limitations}

Despite our best efforts, the presented study is subject to several major limitations. Moderate sample size of 170 participants distributed across 13 of 31 HEIs (participation around $40\%$ of all universities which have mathematics as one of their study cycles) is one major restriction. Another one would be the selection nature of the survey -- while the distribution channel (university Deans) allowed for secure distribution of the questionnaire, it was voluntary in nature. This aspect could have made the resulting answers distribution more specific, as usually people most willing to participate in such studies are the ones with strong opinions on the study topic. Another important limitation was moderate demographic block fill-out (slightly above $50\%$) -- while it was important for us to allow the students to freely express their concerns without having any fear of being identified, it also made the study slightly more difficult to analyse. 
A single cohort subject to the data collection procedure is another issue which contributes to the study limitations.

We also want to point out that we are aware of the reliance on declared intentions rather than realized behaviour which clearly stands out in the conclusive part of the questionnaire.\footnote{Although discrepancies between self-reported behaviour and observed actions are subject to extensive studies (see, e.g. \cite{paschalwebb,paulconner,takeuchi2024}), the meta-analysis of Gnambs and Kaspar \cite{Gnambs2014} shows that people are a lot more eager to report sensitive behaviour via computerized surveys.}
While it might be perceived as a limitation, we believe that tracking the individual anonymous respondents would not be possible without leaving a direct trace leading to them (especially in the case where they've decided to change the university or give up on the second cycle of studies at all, as a way of contacting them after the initial survey would be needed). This in turn could severely impact honesty of the negative viewpoints expressed in the questionnaire, which was clearly a phenomenon we wanted to avoid. 

While cross-sectional nature of the presented study allows for identification of the interdependencies in the data, a longitudinal research in this spirit could possibly allow for identifying the possible shifts in the students' views on learning academic-level mathematics and the market value of a respective degree.

Lastly, at the time of the questionnaire distribution the AI technology was not as prominent and ever-present as it is currently. The changes in education related to the AI usage are beginning to emerge and they tend to escalate rapidly. Studies dedicated to investigation of this phenomenon have started to appear recently and involve various areas of mathematical education, ranging from AI-tutors (Intelligent Tutoring Systems and adaptive platforms), grading and proving feedback, Problem-Based Learning and general problem solving to even curriculum planning.\footnote{See \cite{Torres2024Updating} as well as surveys \cite{Awang2025slr,AIBased_2024} and other paper which discuss further possible directions of research and applications of AI in mathematical education \cite{ssm_walk,Tang_2025,Vanks_2024}.}
This research does not engage this topic in explicit manner, although we acknowledge that rapidly developing AI landscape will most likely enforce subsequent changes into the shifting sands of higher mathematical education. The question to what extent these technological changes will affect our teaching regimes remains at least partially open for now.

\section{Discussion}

Across contexts, mathematics and STEM-students drop-out seems to be driven by the combination of weak prior preparation, early failure in theoretically-demanding math courses, low or worsening interest and self-confidence, anxiety and alienation toward university mathematics, and unsupportive study conditions. Students who persist tend to have stronger foundations, higher interest and self-concept, better learning management, and access to effective institutional support. There are but a few studies which tackle the problem directly, but nevertheless we would like to engage in the discussion with authors of similar studies, even if placed in slightly different regime than ours. Our findings seem to align with their findings, while still introducing several new insights. 

The drop-out factors identified in \cite{heublein2007ursachen}
and discussed in \cite{Geisler2018} (both based on Germany experiences with mathematics drop-out phenomena) focus (among others) on motivation aspects. In \cite{Geisler2018} the author points out that most impactful elements which increase the drop-out probability are the ones related to the motivation (i.e., internal motivations increase the likelihood of staying on the mathematical track, the external ones decrease it).

As far as the drop-out reasons are concerned, \cite{heublein2007ursachen, Geisler2018} both point out  that the ones appearing most are lack of motivation and  work overload, jointly responsible for over half of the drop-outs analysed in the scope of this study. This stays partially consistent with observations from Table \ref{table:reasons} with some important discrepancies and two key differences. The discussed papers strictly referred to the main drop-out reason pointed by the students, while our study asks for main reasons behind the decision to discontinue second-cycle studies on their original university. But the first major difference which strikes out is the fact that the response "I see no point in continuing studying mathematics from the perspective of my professional career." was the most commonly picked amongst all participants in Poland, while German studies had only a few percent of respondents who represent such statement. Second one is the fact that organizational and technical aspects of studies were commonly mentioned by Polish students, at the same time contributing only to several percent (13\%) of the reasons pointed out by their German colleagues. There is a rational explanation for this phenomenon and its roots are in the differences in studied populations. The referenced framework most likely includes broader range of higher-education populations, whereas our data targets specifically mathematics students on the verge of continuing their second-cycle studies. These groups differ systematically in two aspects:
\begin{itemize}
\item The students who are main target of this study have already passed the initial filtering stage of their higher education, where financial, health, and adjustment problems dominate.
\item The population researched in this paper consists of older and therefore more career-oriented people -- which shifts the decision-making calculations toward professional, individual utility rather than adaptation difficulties.
\end{itemize}
Our findings reported in Table \ref{table:reasons} support the points raised by Burke \cite{burke2019student}, who argued that one of the main drivers of student retention is their experience with interactions within their academic and social systems. The questionnaire results indicate that students who choose to change their university or abandon mathematical studies at all are usually reporting low ratio of engagement necessary for such studies to benefits they grant. It is important to notice that the analyses from \cite{burke2019student} were dedicated to students in general higher education, not specifically to mathematical studies.

A difficult conclusion emerging from the presented analyses is that if students report that their standards of knowledge and organization has not been met by the university, we cannot just turn a blind eye. Pretending that the second-cycle studies are for the best and the drop-outs are just the ones who failed to adjust to academic environment and culture is the approach which most likely will cause further problems with overall number of students striving for mathematical degrees which, in the long run, could negatively affect entire STEM-related educational ecosystem in Poland.\footnote{Which is already severely lacking qualified mathematics teachers \cite{Kaminska2023TEACHERS}. By August 2025 Press Agency of Związek Nauczycielstwa Polskiego (Teachers Association in Poland) postulated that there are around 20000 of teachers vacancies which tend to focus around bigger cities.} 

\section{Concluding remarks and recommendations}\label{sec:recommend}

After reviewing and analysing the results of the questionnaire, we believe that we have arrived at certain conclusions, which may improve the current status of second-cycle studies in Poland. Therefore, we advocate for the following set of recommendations:

\begin{itemize}
\item A lot of students tend to be perplexed whether the courses they are attending actually give them any kind of hard skills which would be valuable from the viewpoint of a prospective employer. This is also the main reason for the reluctance of undergraduates towards the second-cycle studies in mathematics. This shows that the universities should place more stress on the cooperation with the business world outside academia and this kind of cooperation should be promoted by the Ministry of Education. Usually the mathematical studies are not expected to provide their graduates with a complete set of skills matching the occupational outlooks. Nevertheless, mathematics is often presented to tentative students as a field of study which allows to flourish on job market -- this prospect should be at least partially matched with a selection of possibly elective, practical courses equipping the students with proper qualifications.
\item There is a component of students of mathematics who believe that the general courses have tendency to be too theoretical. This impression might be an outcome of misplaced expectations after the secondary-level education of mathematics -- after all, the academic education of mathematics differs a lot from the way it is taught at high-schools.
\item Algorithm of allocating public funds to HEIs in Poland \cite{Law2018,Regulation2019} considers multiple criteria, including the number of students, number of staff members and the student-to-staff ratio, research activities undertaken by the universities, doctoral students etc. Amongst those criteria there is no explicit impact of the quality of didactics on the fund allocation. Although the teaching standards do affect the number of applicants for the studies in the long run, there is no direct influence of educational excellence on the funding algorithm. This results in a lot of students being dissatisfied with the substantive and organizational aspects of their stay at their university. Despite the fact that placing additional stress on the quality of teaching might improve the number of students over time, it might be beneficial from the public perspective to introduce some additional markers for education quality into the scope of funds-allocating algorithm.
\item From the results of the survey we can hypothesise that a part of students unsatisfied with their studies might have had different expectations regarding the mathematical higher education. These presumptions might come from completely dissimilar approach to mathematics which is taught at their previous stages of education. One possible solution to this problem is giving the prospective students more insight on the nature of higher mathematics. This might decrease the number of students at the first-cycle studies, but at the same time it will decrease the \textit{shock} factor of sudden switch in their approach to mathematics. This can possibly result in more students who actually want to focus on deepening their knowledge in this field of study and, consequently, higher retention rate between first and second cycle of studies.
\item Lastly, in relationship to the first point of our recommendations, a systematic alumni involvement in curriculum consultation processes is advisable. Academic teachers who have limited connection to the external job market tend to ignore or underestimate the changes which take place in the employers' expectations. A regular external check-up with the graduates who work in the industry may increase transparency of potential career paths of mathematics students and allow them to align their expectations with the reality of the market. Such collaboration can also broaden the promotion of internships and project-based collaboration with external institutions and companies, leading to improved initial experience for students finding their place in the world outside of the academia.
\end{itemize}

The obtained results suggest that the contemporary retention problem in mathematics cannot be reduced solely to academic difficulty or insufficient preparation of students. Instead, the phenomenon appears to emerge from a complex interaction between labour-market expectations, perceived usefulness of studies, organizational aspects of university functioning, specialization availability and students' evolving professional identities. In particular, the discrepancy between appreciation for lecturers' preparation and dissatisfaction with the practical utility of studies indicates that the problem may not necessarily lie in teaching quality itself, but rather in the perceived disconnect between academic mathematics and professional trajectories available to graduates.

An especially noteworthy observation is the partial tension between subjective student perceptions and objective labour-market indicators obtained from the ELA system. While many respondents perceive limited economic justification for continuing second-cycle studies, the administrative data simultaneously demonstrate that the Master's degree still provides measurable labour-market advantages in certain dimensions. This discrepancy itself may constitute an important educational and communication challenge for universities.

At the same time, it should be emphasized that the aim of this paper is not to argue against the Bologna model itself, nor to idealize the former unified system of higher education. Rather, our intention is to identify several structural tensions that emerged within the contemporary educational environment and to contribute to the ongoing discussion on how mathematical education can better respond to both academic and societal expectations.

If we expect mathematics to remain one of the fundamental pillars of technological and scientific development, then maintaining students' long-term engagement with advanced mathematical education should become not merely an institutional concern of individual universities, but an important strategic objective for the broader higher-education system.

\subsection{Further research}

The present study opens several possible directions for subsequent research, both methodological and substantive in nature.

First of all, a natural continuation of this work would involve conducting a longitudinal study tracking the same cohorts of students throughout the entire educational process. While the present questionnaire captures declared intentions regarding continuation of studies, linking these declarations with realized educational trajectories would allow for a more precise assessment of the discrepancy between intentions and actual behaviour. Such an approach could also reveal the temporal dynamics of motivation, academic self-concept and perceived usefulness of mathematical education.

Secondly, the strong inter-university variability visible in the ELA-based analyses suggests that retention in mathematics is at least partially an institution-level phenomenon. Future studies could therefore focus on identifying organizational and curricular factors characteristic of universities with particularly high retention rates. This could include comparisons of: curriculum structures; ratios of theoretical to practical courses; available specializations; forms of cooperation with industry;
 thesis supervision systems; organizational flexibility for working students; and support mechanisms for students.

Another important direction concerns the relationship between mathematical studies and labour-market expectations. The current results suggest that many students perceive insufficient professional utility in continuing second-cycle mathematical education. Future research in this area could therefore investigate which  competencies employers actually expect from mathematics graduates and to what extent university curricula reflect these expectations. A separate question in this spirit could be whether a closer integration between universities and industry partners improves retention and graduate satisfaction.

The role of specialization choice also deserves a more detailed investigation. Our findings (cf. Table \ref{table:postbhcorec}) indicate that satisfaction with the chosen specialization is one of the strongest differentiating factors between students willing and unwilling to continue second-cycle studies. Future work could therefore examine how specialization structures influence retention and whether early exposure to specialized courses has any effect on students' motivation levels.

A particularly promising extension of this study could possibly involve qualitative research methods. While quantitative questionnaires allow for detecting large-scale patterns, in-depth interviews or focus-group studies could provide a substantially richer understanding of students' perceptions of academic mathematics as well as their individual expectations regarding university education.

Finally, the rapidly evolving influence of artificial intelligence on higher education constitutes another important research direction one cannot overlook. Since the data collection for this study preceded the current large-scale proliferation of generative AI systems, the presented results should be interpreted cautiously within a rapidly changing educational landscape. Future studies could investigate the impact of AI tools availability on students' mathematical self-concept and critical thinking skills.

\subsection*{Acknowledgements}

Firstly, we would like to extend our utmost gratitude to all of the people who enabled us to conduct this research -- Deans and Directors of multiple mathematically-oriented institutes and faculties of universities in Poland. Without your precious assistance this article would not see the light of the day.

Secondly, we would like to thank once again all the students participating in this survey, especially those who shared additional details on their stories, decisions they made and the situations they have faced during their studies. This paper is but a speck of your story, but hopefully it will make an impact on higher education in Poland.

We would also like to thank all the people, belonging to neither of the previous groups, whose suggestions led to improved quality of this paper. This includes Piotr Nowakowski, Szymon Głąb, Marta Borowska-Stefańska, Dominik Pawlak and many others.

Finally, Filip Turoboś would like to apologize for the delay caused in publishing of the final version of this paper and express his utmost thanks for all the support he received from his co-authors, family and friends during his cancer treatment. 

\begin{center}
\textit{Two are better than one, because they have a good return for their labor. If either of them falls down, one can help the other up. But pity anyone who falls and has no one to help them up.}
\end{center}

{\hspace{0.7\textwidth}{Ecclesiastes 4:9-10 NIV}}

\subsection*{Contributions}

Conceptualization and methodology: all the authors; software: F.T.; validation: Ż.F.; formal analysis: all the authors; investigation: F.T.; preliminary data collection: J.S., N.M. and F.T.; main data curation: F.T.;statistical verification: F.T. and Ż.F.; writing -- original draft preparation: F.T.; writing -- review and editing: all the authors; visualization, F.T.; supervision: F.T.; project administration: F.T.; funding acquisition: F.T. 


\subsection*{Funding}
F. Turoboś research on this topic was funded by the National Science Centre, Poland, under the program 
\textbf{MINIATURA 7}, grant number: 2023/07/X/ST1/01687. 

\subsection*{Data availability statement}

The cleaned and anonymized dataset generated and analysed during the current study is publicly available in a GitHub repository at (\href{https://github.com/janfifian/SecondCycleStudiesDataset}{https://github.com/janfifian/SecondCycleStudiesDataset}) and archived on Zenodo at \href{https://doi.org/10.5281/zenodo.17692881}{https://doi.org/10.5281/zenodo.17692881} under the Creative Commons Attribution 4.0 International license.

The ELA graduate-tracking dataset utilized in this study is publicly accessible via the Polish Graduate Tracking System website:
 \href{https://www.ela.nauka.gov.pl/en}{https://www.ela.nauka.gov.pl/en}.

\bibliography{bibliografia}

@article{OlenczukPaszel2009,
  author       = {Anna Oleńczuk-Paszel},
  title        = {Realizacja Procesu Bolońskiego w Polsce a rozwój kapitału społecznego},
  journal      = {Folia Pomeranae Universitatis Technologiae Stetinensis. Oeconomica},
  volume       = {273},
  number       = {56},
  pages        = {155--162},
  year         = {2009},
  note         = {(available in Polish)}
}

@misc{PrawoOSzkolnictwieWyzszym2005,
  title        = {Ustawa z 27 lipca 2005 r. Prawo o szkolnictwie wyższym},
  year         = {2005},
  note         = {DzU z 2005 r., nr 164, poz. 1365 z późn. zm.},
  url          = {https://eli.gov.pl/eli/DU/2005/1365/ogl}
}

@misc{SorbonneDeclaration1998,
  title        = {Sorbonne Joint Declaration: Joint declaration on harmonisation of the architecture of the European higher education system},
  year         = {1998},
  note         = {Signed in Paris, May 25, 1998},
  url          = {https://ehea.info/media.ehea.info/file/1998_Sorbonne/61/2/1998_Sorbonne_Declaration_English_552612.pdf}
}

@misc{BolognaDeclaration1999,
  title        = {The Bologna Declaration: Joint Declaration of the European Ministers of Education},
  year         = {1999},
  note         = {Signed in Bologna, June 19, 1999},
  url          = {https://ehea.info/media.ehea.info/file/Ministerial_conferences/02/8/1999_Bologna_Declaration_English_553028.pdf}
}

@inproceedings{kilian2020predicting,
  title={Predicting math student success in the initial phase of college with sparse information using approaches from statistical learning},
  author={Kilian, Pascal and Loose, Frank and Kelava, Augustin},
  booktitle={Frontiers in Education},
  volume={5},
  pages={502698},
  year={2020},
  organization={Frontiers Media SA},
  doi={10.3389/feduc.2020.502698}
}

@article{burke2019student,
  title={Student retention models in higher education: A literature review},
  author={Burke, Adam},
  journal={College and University},
  volume={94},
  number={2},
  pages={12--21},
  year={2019},
  publisher={American Association of Collegiate Registrars and Admissions Officers}
}

@article{tight2020student,
  title={Student retention and engagement in higher education},
  author={Tight, Malcolm},
  journal={Journal of further and Higher Education},
  volume={44},
  number={5},
  pages={689--704},
  year={2020},
  publisher={Taylor \& Francis},
  doi={https://doi.org/10.1080/0309877X.2019.1576860}
}

@article{xavier2020literature,
  title={A literature review on the definitions of dropout in online higher education},
  author={Xavier, Marlon and Meneses, Julio},
  year={2020},
  journal={European Distance and E-Learning Network (EDEN)},
  doi={http://doi.org/10.38069/edenconf-2020-ac0004},
  pages={73-80}
}

@article{valencia2024dropout,
  title={Dropout in postgraduate programs: a underexplored phenomenon -- a scoping review},
  author={Valencia Quecano, Lira Isis and Guzm{\'a}n Rinc{\'o}n, Alfredo and Barrag{\'a}n Moreno, Sandra},
  journal={Cogent Education},
  volume={11},
  number={1},
  pages={2326705},
  year={2024},
  publisher={Taylor \& Francis},
  doi={https://doi.org/10.1080/2331186X.2024.2326705}
}

@article{Mller2022,
  title = {Social Inequality in Dropout from Higher Education in Germany. Towards Combining the Student Integration Model and Rational Choice Theory},
  volume = {64},
  ISSN = {1573-188X},
  url = {http://dx.doi.org/10.1007/s11162-022-09703-w},
  DOI = {10.1007/s11162-022-09703-w},
  number = {2},
  journal = {Research in Higher Education},
  publisher = {Springer Science and Business Media LLC},
  author = {M\"{u}ller,  Lars and Klein,  Daniel},
  year = {2022},
  pages = {300–330}
}

@article{nieuwoudt2021whystudents,
  author = {Nieuwoudt, Johanna E. and Pedler, Megan L.},
  title = {Student Retention in Higher Education: Why Students Choose to Remain at University},
  journal = {Journal of College Student Retention: Research, Theory \& Practice},
  year = {2021},
  pages = {326-349},
  number = {2},
  volume = {25},
  doi = {10.1177/1521025120985228}
}

@article{IJESL2237,
	author = {Othman Aljohani},
	title = {A Review of the Contemporary International Literature on Student Retention in Higher Education},
	journal = {International Journal of Education and Literacy Studies},
	volume = {4},
	number = {1},
	year = {2016},
	pages = {40--52},	
	url = {https://journals.aiac.org.au/index.php/IJELS/article/view/2237},
	doi = {10.7575/aiac.ijels.v.4n.1p.40}
}

@article{fenwick2009recruitment,
  title={Recruitment and retention of mathematics students in Canadian universities},
  author={Fenwick-Sehl, Laura and Fioroni, Marcella and Lovric, Miroslav},
  journal={International Journal of Mathematical Education in Science and Technology},
  volume={40},
  number={1},
  pages={27--41},
  year={2009},
  publisher={Taylor \& Francis},
  doi={10.1080/00207390802568192}
}

@article{Lytle2023STEMEngagement,
  author    = {Ashley Lytle and Jiyun Elizabeth L. Shin},
  title     = {Self and Professors’ Incremental Beliefs as Predictors of STEM Engagement Among Undergraduate Students},
  journal   = {International Journal of Science and Mathematics Education},
  year      = {2023},
  volume    = {21},
  number    = {3},
  pages     = {1013--1029},
  doi       = {10.1007/s10763-022-10272-8},
  url       = {https://doi.org/10.1007/s10763-022-10272-8}
}

@article{rotem2021dropping,
  title={Dropping out of master’s degrees: Objective predictors and subjective reasons},
  author={Rotem, Nir and Yair, Gad and Shustak, Elad},
  journal={Higher Education Research \& Development},
  volume={40},
  number={5},
  pages={1070--1084},
  year={2021},
  publisher={Taylor \& Francis},
  doi={https://doi.org/10.1080/07294360.2020.1799951}
}

@article{         harris2020array,
 title         = {Array programming with {NumPy}},
 author        = {Charles R. Harris and K. Jarrod Millman and et al.},
 year          = {2020},
 month         = sep,
 journal       = {Nature},
 volume        = {585},
 number        = {7825},
 pages         = {357--362},
 doi           = {10.1038/s41586-020-2649-2},
 publisher     = {Springer Science and Business Media {LLC}},
 url           = {https://doi.org/10.1038/s41586-020-2649-2}
}

@inproceedings{pandas,
  author    = { {W}es {M}c{K}inney },
  title     = { {D}ata {S}tructures for {S}tatistical {C}omputing in {P}ython },
  booktitle = { {P}roceedings of the 9th {P}ython in {S}cience {C}onference },
  pages     = { 56 - 61 },
  year      = { 2010 },
  editor    = { {S}t\'efan van der {W}alt and {J}arrod {M}illman },
  doi       = { 10.25080/Majora-92bf1922-00a }
}

@article{scipy,
       author = {{Virtanen}, Pauli and {Gommers}, Ralf and et al.},
        title = "{SciPy 1.0: Fundamental Algorithms for Scientific
                  Computing in Python}",
      journal = {Nature Methods},
      year = "2020",
      volume={17},
      pages={261--272},
      adsurl = {https://rdcu.be/b08Wh},
      doi = {https://doi.org/10.1038/s41592-019-0686-2},
}

@conference{jupyter,
  title={Jupyter Notebooks-a publishing format for reproducible computational workflows},
  author={KLUYVER{\textordfeminine}$^1$, Thomas and RAGAN-KELLEYb$^1$, Benjamin and P{\'e}rez, Fernando and Granger, Brian and Bussonnier, Matthias and Frederic, Jonathan and Kelley, Kyle},
  booktitle={Positioning and power in academic publishing: players, agents and agendas: proceedings of the 20th International Conference on Electronic Publishing},
  pages={87},
  year={2016},
  organization={IOS press},
  doi={https://doi.org/10.3233/978-1-61499-649-1-87}
}

@article{matplotlib,
  author={J. D. {Hunter}},
  journal={{C}omputing in {S}cience \& {E}ngineering}, 
  title={{Matplotlib: A 2D Graphics Environment}}, 
  year={2007},
  volume={9},
  number={3},
  pages={90-95}}

@article{Waskom2021,
    doi = {10.21105/joss.03021},
    url = {https://doi.org/10.21105/joss.03021},
    year = {2021},
    publisher = {The Open Journal},
    volume = {6},
    number = {60},
    pages = {3021},
    author = {Michael L. Waskom},
    title = {seaborn: statistical data visualization},
    journal = {Journal of Open Source Software}
 }

@online{plotly, author = {Plotly Technologies Inc.}, title = {Collaborative data science}, publisher = {Plotly Technologies Inc.}, address = {Montreal, QC}, year = {2015}, url = {https://plot.ly} }

@article{kaplan2020multi,
  title={Multi-mode question pretesting: Using traditional cognitive interviews and online testing as complementary methods},
  author={Kaplan, Robin L and Edgar, Jennifer},
  journal={Survey Methods: Insights from the Field},
  pages={1--14},
  year={2020},
  publisher={DEU}
}

@article{hilton2017importance,
  title={The importance of pretesting questionnaires: a field research example of cognitive pretesting the Exercise referral Quality of Life Scale (ER-QLS)},
  author={Hilton, Charlotte Emma},
  journal={International Journal of Social Research Methodology},
  volume={20},
  number={1},
  pages={21--34},
  year={2017},
  publisher={Taylor \& Francis}
}

@article{neuert2021effects,
  title={Effects of the number of open-ended probing questions on response quality in cognitive online pretests},
  author={Neuert, Cornelia E and Lenzner, Timo},
  journal={Social Science Computer Review},
  volume={39},
  number={3},
  pages={456--468},
  year={2021},
  publisher={SAGE Publications Sage CA: Los Angeles, CA}
}

@article{rothgeb2007questionnaire,
  title={Questionnaire pretesting methods: Do different techniques and different organizations produce similar results?},
  author={Rothgeb, Jennifer and Willis, Gordon and Forsyth, Barbara},
  journal={Bulletin of Sociological Methodology/Bulletin de M{\'e}thodologie Sociologique},
  volume={96},
  number={1},
  pages={5--31},
  year={2007},
  publisher={Sage Publications Sage CA: Thousand Oaks, CA}
}

@article{kaczmirek2017higher,
  title={Higher data quality in web probing with EvalAnswer: a tool for identifying and reducing nonresponse in openended questions},
  author={Kaczmirek, Lars and Meitinger, Katharina and Behr, Doroth{\'e}e},
  year={2017},
  publisher={DEU},
  journal = {GESIS Papers},
  volume = {2017/01},
  doi = {https://doi.org/10.21241/ssoar.51100}
}

@article{geisen2020compendium,
  title={A compendium of web and mobile survey pretesting methods},
  author={Geisen, Emily and Murphy, Joe},
  journal={Advances in questionnaire design, development, evaluation and testing},
  pages={287--314},
  year={2020},
  publisher={Wiley Online Library}
}

@inproceedings{du2016research,
  title={Research and improvement of Apriori algorithm},
  author={Du, Jiaoling and Zhang, Xiangli and Zhang, Hongmei and Chen, Lei},
  booktitle={2016 Sixth International Conference on Information Science and Technology (ICIST)},
  pages={117--121},
  year={2016},
  organization={IEEE}
}

@inproceedings{agrawal1993mining,
  title={Mining association rules between sets of items in large databases},
  author={Agrawal, Rakesh and Imieli{\'n}ski, Tomasz and Swami, Arun},
  booktitle={Proceedings of the 1993 ACM SIGMOD international conference on Management of data},
  pages={207--216},
  year={1993}
}

@book{hill2010so,
  title={Why so few? Women in science, technology, engineering, and mathematics.},
  author={Hill, Catherine and Corbett, Christianne and St Rose, Andresse},
  year={2010},
  publisher={ERIC}
}

@article{silbey2016so,
  title={Why do so many women who study engineering leave the field},
  author={Silbey, Susan S.},
  journal={Harvard Business Review},
  volume={23},
  pages={1--8},
  year={2016},
  url={https://hbr.org/2016/08/why-do-so-many-women-who-study-engineering-leave-the-field}
}

@article{Kaminska2023TEACHERS,
title={Teachers’ working conditions in times of professional crisis: evidence of Poland},
author={Małgorzata Kamińska},
journal={The Modern Higher Education Review},
year={2023},
doi={10.28925/2617-5266.2023.83}
}

@article{espinoza2023whydo,
  author = {Espinoza, Oscar and González, Luis Eduardo and Sandoval, Luis and McGinn, Noel and Corradi, Bruno},
  title = {Why do students leave? Persistence in selective universities},
  journal = {Research Papers in Education},
  year = {2023},
  doi = {10.1080/02671522.2023.2238274}
}

@book{Crosier2013,
  author    = {David Crosier and Teodora Parveva},
  title     = {The Bologna Process: Its Impact in Europe and Beyond},
  year      = {2013},
  publisher = {UNESCO International Institute for Educational Planning},
  address   = {Paris},
  series    = {Fundamentals of Educational Planning},
  volume    = {97},
  isbn      = {978-92-803-1368-0},
  url       = {https://unesdoc.unesco.org/ark:/48223/pf0000220649}
}

@article{Terry2008,
  author    = {Laurel S. Terry},
  title     = {The Bologna Process and its Impact in Europe: It's So Much More Than Degree Changes},
  journal   = {Vanderbilt Journal of Transnational Law},
  volume    = {41},
  pages     = {107},
  year      = {2008},
  note      = {Available at SSRN: \url{https://ssrn.com/abstract=1139805}}
}

@incollection{Huisman2012,
  author    = {Huisman, J. and Adelman, C. and Hsieh, C.C. and Shams, F. and Wilkins, S.},
  title     = {Europe’s Bologna process and its impact on global higher education},
  booktitle = {Handbook of International Higher Education},
  editor    = {Deardorff, D.K. and de Wit, H. and Heyl, J.D. and Adams, T.},
  pages     = {81--100},
  year      = {2012},
  publisher = {Sage Publications},
  address   = {Thousand Oaks}
}

@inproceedings{Usaci2012,
  author    = {Doina Usaci and Mariana Norel and Rodica Mariana Niculescu and Daciana Lupu},
  title     = {Implication of Bologna Process on the Academic Curriculum: An Approach Focused on Student’s Perspective},
  booktitle = {Handbook of International Higher Education},
  pages     = {81--100},
  year      = {2012},
  publisher = {Sage Publications},
  address   = {Thousand Oaks}
}

@misc{Law2018,
  title     = {Ustawa z dnia 20 lipca 2018 r. - Prawo o szkolnictwie wyższym i nauce},
  year      = {2018},
  url       = {https://isap.sejm.gov.pl/isap.nsf/DocDetails.xsp?id=WDU20180001668}
}

@misc{Regulation2019,
  title     = {Rozporządzenie Ministra Nauki i Szkolnictwa Wyższego z dnia 9 września 2019 r. w sprawie sposobu podziału środków finansowych na utrzymanie i rozwój potencjału dydaktycznego oraz potencjału badawczego znajdujących się w dyspozycji ministra właściwego do spraw szkolnictwa wyższego i nauki oraz na zadania związane z utrzymaniem powietrznych statków szkolnych i specjalistycznych ośrodków szkoleniowych kadr powietrznych},
  year      = {2019},
  url       = {https://isap.sejm.gov.pl/isap.nsf/DocDetails.xsp?id=WDU20190001777}
}

@article{Torres2024Updating,
title={Updating Calculus Teaching with AI: A Classroom Experience},
author={R. Torres-Peña and Darwin Peña-González and Ellery Chacuto-López and E. A. Ariza and Diego Vergara},
journal={Education Sciences},
year={2024},
doi={10.3390/educsci14091019}
}

@article{Awang2025slr,
  title     = {Current practices and future direction of artificial
               intelligence in mathematics education: A systematic review},
  author    = {Awang, Liz A. and Yusop, Farrah D. and Danaee, Mahmoud},
  journal   = {Int. Electron. J. Math. Educ.},
  publisher = {Modestum Ltd},
  volume    =  {20},
  number    =  {2},
  pages     = {em0823},
  doi = {https://doi.org/10.29333/iejme/16006},
  year      =  {2025}
}

@article{AIBased_2024, 
title = {AI-Based Tools in Mathematics Education: A Systematic Review of Characteristics, Applications, and Evaluation Methods},
author = {KP Mredula and Roman Jonita and Priti Sajja},
volume={2}, 
url={https://irjaeh.com/index.php/journal/article/view/308}, 
DOI={10.47392/IRJAEH.2024.0268}, 
number={07}, 
journal={International Research Journal on Advanced Engineering Hub (IRJAEH)}, 
year={2024}, 
pages={1958--1967} 
}

@article{Tang_2025, 
title={Artificial Intelligence in Mathematics Education: Trends, Challenges, and Opportunities}, 
volume={3}, 
url={https://ejournal.uinsaizu.ac.id/index.php/ijrme/article/view/13496}, 
DOI={10.24090/ijrme.v3i1.13496}, 
number={1}, 
journal={International Journal of International Journal of Science and Mathematics Education}, 
author={Tang, William Ko-Wai}, 
year={2025}, 
month={Jul.}, 
pages={75--90} 
}

@article{ssm_walk,
author = {Walkington, Candace},
title = {The implications of generative artificial intelligence for mathematics education},
journal = {School Science and Mathematics},
doi = {https://doi.org/10.1111/ssm.18356},
url = {https://onlinelibrary.wiley.com/doi/abs/10.1111/ssm.18356},
eprint = {https://onlinelibrary.wiley.com/doi/pdf/10.1111/ssm.18356},
pages={1--10},
year={2025}
}

@article{Vanks_2024, 
title={Generative Artificial Intelligence on Mobile Devices in the University Preparation of Future Teachers of Mathematics}, 
volume={18}, 
url={https://online-journals.org/index.php/i-jim/article/view/51221}, 
DOI={10.3991/ijim.v18i18.51221}, 
number={18}, 
journal={International Journal of Interactive Mobile Technologies (iJIM)}, 
author={Vankúš, Peter}, 
year={2024}, 
pages={19--33} }

@article{Gnambs2014,
  title = {Disclosure of sensitive behaviors across self-administered survey modes: a meta-analysis},
  volume = {47},
  ISSN = {1554-3528},
  url = {http://dx.doi.org/10.3758/s13428-014-0533-4},
  DOI = {10.3758/s13428-014-0533-4},
  number = {4},
  journal = {Behavior Research Methods},
  publisher = {Springer Science and Business Media LLC},
  author = {Gnambs,  Timo and Kaspar,  Kai},
  year = {2015},
  pages = {1237–-1259}
}

@article{takeuchi2024,
author = {Makito Takeuchi and Junichiro Niimi and Takahiro Hoshino},
title ={Handling the Inconsistency between Self-Report and the Actual Behavior: Validity of Excluding Survey Participants with Insufficient Effort Responding},
journal = {International Journal of Market Research},
volume = {66},
number = {4},
pages = {451--472},
year = {2024},
doi = {10.1177/14707853231209933},
URL = {https://doi.org/10.1177/14707853231209933},
eprint = {https://doi.org/10.1177/14707853231209933}
}

@article{paschalwebb,
author = {Sheeran, Paschal and Webb, Thomas L.},
title = {The Intention–Behavior Gap},
journal = {Social and Personality Psychology Compass},
volume = {10},
number = {9},
pages = {503-518},
doi = {https://doi.org/10.1111/spc3.12265},
url = {https://compass.onlinelibrary.wiley.com/doi/abs/10.1111/spc3.12265},
eprint = {https://compass.onlinelibrary.wiley.com/doi/pdf/10.1111/spc3.12265},
year = {2016}
}

@ARTICLE{paulconner,
AUTHOR={Conner, Mark  and Norman, Paul },
TITLE={Understanding the intention-behavior gap: The role of intention strength},
JOURNAL={Frontiers in Psychology},
VOLUME={Volume 13 - 2022},
YEAR={2022},
URL={https://www.frontiersin.org/journals/psychology/articles/10.3389/fpsyg.2022.923464},
DOI={10.3389/fpsyg.2022.923464},  
ISSN={1664-1078}
}

@article{Benjamini1995Controlling,
title={Controlling the false discovery rate: a practical and powerful approach to multiple testing},
author={Y. Benjamini and Y. Hochberg},
journal={Journal of the royal statistical society series b-methodological},
year={1995},
volume={57},
pages={289-300},
doi={10.1111/j.2517-6161.1995.tb02031.x}
}

@article{Barber2017The,
title={The p-filter: multilayer false discovery rate control for grouped hypotheses},
author={R. Barber and Aaditya Ramdas},
journal={Journal of the Royal Statistical Society: Series B (Statistical Methodology)},
year={2017},
volume={79},
doi={10.1111/rssb.12218}
}

@inproceedings{rach2019self,
  title={Self-concept in university mathematics courses},
  author={Rach, Stefanie and Ufer, Stefan and Kosiol, Timo},
  booktitle={Eleventh Congress of the European Society for Research in Mathematics Education},
  number={18},
  year={2019},
  organization={Freudenthal Group; Freudenthal Institute; ERME},
  url={https://hal.science/hal-02410205v1}
}

@article{Geisler2023,
  title = {The relation between attitudes towards mathematics and dropout from university mathematics—the mediating role of satisfaction and achievement},
  volume = {112},
  ISSN = {1573-0816},
  url = {http://dx.doi.org/10.1007/s10649-022-10198-6},
  DOI = {10.1007/s10649-022-10198-6},
  number = {2},
  journal = {Educational Studies in Mathematics},
  publisher = {Springer Science and Business Media LLC},
  author = {Geisler,  Sebastian and Rach,  Stefanie and Rolka,  Katrin},
  year = {2023},
  month = jan,
  pages = {359–381}
}

@incollection{Geisler2018,
author={Geisler, Sebastian},
editor={Rott, Benjamin
and T{\"o}rner, G{\"u}nter
and Peters-Dasdemir, Joyce
and M{\"o}ller, Anne
and Safrudiannur},
title={It's All About Motivation?---A Case Study Concerning Dropout and Persistence in University Mathematics},
bookTitle={Views and Beliefs in Mathematics Education: The Role of Beliefs in the Classroom},
year={2018},
publisher={Springer International Publishing},
address={Cham},
pages={115--124},
isbn={978-3-030-01273-1},
doi={10.1007/978-3-030-01273-1_11},
url={https://doi.org/10.1007/978-3-030-01273-1_11}
}

@article{Almukhambetova2023RetentionSTEM,
  author    = {Ainur Almukhambetova and Aliya Kuzhabekova and Daniel Hernández-Torrano},
  title     = {Hidden Bias, Low Expectations, and Social Stereotypes: Understanding Female Students’ Retention in Math-Intensive STEM Fields},
  journal   = {International Journal of Science and Mathematics Education},
  year      = {2023},
  volume    = {21},
  number    = {2},
  pages     = {535--557},
  doi       = {10.1007/s10763-022-10256-8},
  url       = {https://doi.org/10.1007/s10763-022-10256-8}
}

@article{Almukhambetova2023LeakyPipeline,
  author    = {Ainur Almukhambetova and Daniel Hernández-Torrano and Alexandra Nam},
  title     = {Fixing the Leaky Pipeline for Talented Women in STEM},
  journal   = {International Journal of Science and Mathematics Education},
  year      = {2023},
  volume    = {21},
  pages     = {305--324},
  doi       = {10.1007/s10763-021-10239-1},
  url       = {https://doi.org/10.1007/s10763-021-10239-1}
}

@article{heublein2007ursachen,
  title={Ursachen des Studienabbruchs in Bachelor-und in herk{\"o}mmlichen Studieng{\"a}ngen},
  author={Heublein, Ulrich and Hutzsch, Christopher and Schreiber, Jochen and Sommer, Dieter and Besuch, Georg},
  journal={Ergebnisse einer bundesweiten Befragung von Exmatrikulierten des Studienjahres},
  volume={8},
  number={2},
  year={2009},
  url={https://www.dzhw.eu/pdf/21/studienabbruch_ursachen.pdf}
}

@article{Kocsis2024Factors,
title={Factors influencing academic performance and dropout rates in higher education},
author={Ádám Kocsis and G. Molnár},
journal={Oxford Review of Education},
year={2024},
volume={51},
pages={414 - 432},
doi={10.1080/03054985.2024.2316616}
}

@article{Willson1976Critical,
title={Critical Values of the Rank-Biserial Correlation Coefficient},
author={V. Willson},
journal={Educational and Psychological Measurement},
year={1976},
volume={36},
pages={297 - 300},
doi={10.1177/001316447603600207}
}

@article{Tapio2025The,
title={The Role of Data Assumptions in Selecting Between Parametric and Nonparametric Tests},
author={Rachelle P. Tapio},
journal={Asian Journal of Probability and Statistics},
year={2025},
doi={10.9734/ajpas/2025/v27i11830}
}

@article{Pangesti2025Application,
title={Application of Cluster Analysis and Correlation between Mathematics and Natural Sciences Subject Based on Student Test Scores Using K-Means Clustering},
author={S. Pangesti and Nurul Farisah and Binti Roslan and Andreansyah},
journal={Indonesian Journal of Applied Mathematics and Statistics},
year={2025},
doi={10.71385/idjams.v2i1.8}
}

@article{Anabo2023CORRELATES,
title={CORRELATES OF MATHEMATICS PERFORMANCE OF GRADE 9 LEARNERS IN SECONDARY SCHOOLS DIVISION OF EASTERN SAMAR AMIDST PANDEMIC},
author={Roy O. Anabo},
journal={EPRA International Journal of Multidisciplinary Research (IJMR)},
year={2023},
doi={10.36713/epra13649}
}
\bibliographystyle{plain}

\end{document}